\documentclass[submit]{epsv8}
\usepackage{lineno}
\usepackage{amssymb}
\usepackage{amsmath}
\usepackage{color}
\usepackage{graphicx}
\usepackage{xcolor}

\newcommand {\reviser}[1]{\textcolor{black}{#1}}
\newcommand {\reviseb}[1]{\textcolor{black}{#1}}
\newcommand {\revrr}[1]{\textcolor{black}{#1}}
\newcommand {\vect}[1]{\mbox{\boldmath $#1$}}

\newcommand {\pdif}[3][]{\frac{\partial^{#1}#2}{\partial#3^{#1}}}

\newcommand {\lsim}{\hspace{0.3em}\raisebox{0.4ex}{$<$}\hspace{-0.75em}\raisebox{-.7ex}{$\sim$}\hspace{0.3em}}

\makeatletter
\def\mart{\@ifnextchar[{\mart@@}{\mart@}}
\def\mart@@[#1]#2{\sqrt[#1]{\mathstrut{#2}}}
\def\mart@#1{\sqrt{\mathstrut{#1}}}
\makeatother
\newcommand {\Alfven}{Alfv\'{e}n}

\newcommand{\sgn}{{\rm sgn}}

\newcommand{\median}{{\rm median}}
\newcommand{\Ampere}{{Amp{\`e}re }}

\newcommand{\apj}{Astrophysical Journal }
\newcommand{\apjs}{Astrophysical Journal Supplement Series }
\newcommand{\jgr}{Journal of Geophysical Research }

\title{A high-order weighted positive and flux conservative method for the Vlasov equation}
\author{
Takashi Minoshima (corresponding author), Center for Mathematical Science and Advanced Technology, Japan Agency for Marine-Earth Science and Technology, 3173-25 Syowa-machi, Kanazawaku, Yokohama 236-0001, Japan, minoshim@jamstec.go.jp\\
Yosuke Matsumoto, Institute for Advanced Academic Research, Chiba University, 1-33 Yayoi-cho, Inage-ku, Chiba 263-8522, Japan, ymatumot@chiba-u.jp
}

\abstract{
We present a high-order conservative, positivity-preserving, and non-oscillatory scheme for solving the Vlasov equation.
The scheme attains formal fifth-order accuracy through a convex combination of positive and non-oscillatory polynomials in substencils. 
Nonlinear weights for these polynomials are formulated that assign higher priority to substencils with larger $L^2$ norm to enhance resolution while maintaining positivity and non-oscillatory properties.
 An approximate dispersion relation indicates that the spectral properties of the present scheme outperform those of an underlying fifth-order scheme and even surpass those of a seventh-order scheme in certain wavenumber ranges. 
We apply this scheme to the one-dimensional Vlasov-\Ampere equations and the two-dimensional Vlasov-Maxwell equations, and demonstrate high-resolution simulations with improved conservation of entropy.
}


\begin{document}

\maketitle

\section{Introduction}\label{sec:introduction}
Kinetic plasma simulations are extensively utilized to investigate nonlinear physics in astrophysical, space, and laboratory plasma environments, including instabilities, wave-particle interactions, and particle acceleration.
Two major ways are employed for the kinetic plasma simulation: Particle-In-Cell (PIC) simulation \citep{PIC} and Vlasov simulation \citep{1976JCoPh..22..330C}.
The Vlasov simulation is an Eulerian approach that discretizes the distribution function to solve the Vlasov equation as an advection equation in phase space.
While the Vlasov simulation can address numerical difficulties inherent in PIC simulations and provide high-resolution solutions, its bottleneck is the substantial computational cost to allocate grid points in up to six-dimensional phase space.
As the number of grid points per dimension is severely restricted, typically $\lsim O(10^2)$ on currently available computer resources, high-order schemes are essential for accurately solving the Vlasov equation.

The semi-Lagrangian scheme for the advection equation is utilized in Vlasov simulations \citep{1999JCoPh.149..201S,1999CoPhC.120..122N,2009CoPhC.180.1730C,2011JCoPh.230.6800M,2013JCoPh.236...81M}.
This method updates the profile by advecting it from an upstream position along the characteristic curve, where the upstream profile is reconstructed through interpolation.
\reviser{The scheme is advantageous that its time step is not constrained by the CFL condition once the upstream position is determined.}
A high-order extension of the scheme is attained by using high-degree polynomials for interpolation.
\reviseb{Recently, novel high-order semi-Lagrangian schemes have been developed for the Vlasov-\Ampere and Vlasov-Maxwell systems that guarantee energy conservation \citep{liu2023efficient,2025JCoPh.52913858L}.}
However, linear high-order schemes suffer from numerical oscillations, thereby necessitating the use of nonlinear schemes to limit the slope of profiles.
In Vlasov simulations, nonlinear schemes are implemented to suppress numerical oscillations that can cause unphysical generation of plasma waves and negative densities in phase space. 
A drawback of using nonlinear schemes is inherent numerical diffusion, which can obscure physically meaningful profiles and degrade energy conservation.

A variety of nonlinear semi-Lagrangian schemes have been proposed to satisfy key properties required for the Vlasov equation, such as mass conservation, positivity, and non-oscillatory behavior \citep{2001JCoPh.172..166F,2008EP&S...60..773U,2010JCoPh.229.1927C,qiu2010conservative,2012CoPhC.183.1094U,xiong2014high,cai2016conservative,2017ApJ...849...76T,sirajuddin2019truly}.
In this paper, we present a new high-order semi-Lagrangian scheme that is conservative, positivity-preserving, and non-oscillatory, for solving the Vlasov equation.
The scheme builds upon the third-order positive and flux conservative (PFC) method \citep{2001JCoPh.172..166F} and its variant \citep{2008EP&S...60..773U}, extending it to attain fifth-order accuracy.
This is achieved by a convex combination of low-order polynomials in substencils, inspired by the weighted essentially non-oscillatory (WENO) scheme \citep{1996JCoPh.126..202J}.
Unlike conventional WENO schemes, however, we utilize positive and non-oscillatory polynomials derived from the PFC method.
This approach allows us to assign weights to substencils so as to outperform the underlying linear scheme.
In \S \ref{sec:posit-flux-cons}, we briefly review the PFC method and provide several improvements.
Subsequently, we detail the new scheme in \S \ref{sec:weight-posit-flux}.
In \S \ref{sec:numer-exper}, we conduct numerical experiments to assess the performance of our scheme.
Finally, we conclude the paper in \S \ref{sec:conclusion}.

\section{Positive and flux conservative (PFC) method: Review and improvements}\label{sec:posit-flux-cons}
We examine the one-dimensional advection equation for the profile $f(x,t)$ written in a conservative form,
\begin{eqnarray}
\pdif{f}{t} + \pdif{}{x}\left(v f\right) = 0. \label{eq:1}
\end{eqnarray}
The profile is discretized in the $j$-th cell $\in [x_j-\Delta x/2,x_j+\Delta x/2]$ at the time $t_n = n \Delta t$ as
\begin{eqnarray}
f_j^n = \frac{1}{\Delta x}\int_{x_j-\Delta x/2}^{x_j+\Delta x/2}f(x,t_n)dx,\label{eq:2}
\end{eqnarray}
where $\Delta x$ and $\Delta t$ are the grid width and the time step, respectively.
Equation (\ref{eq:1}) can be updated in a conservative semi-Lagrangian manner,
\begin{eqnarray}
 f_j^{n+1} =  f_j^{n} - \Phi_j + \Phi_{j-\sgn(v)},\label{eq:3}
\end{eqnarray}
where $\Phi_j$ is the outgoing flux to a downstream cell,
\begin{eqnarray}
 \Phi_j &=& \left\{
\begin{array}{l}
 \int_{0.5-\xi}^{0.5} F_j(x) dx, \;\;\; {\rm for} \; v>0, \\
 \int_{-0.5}^{-0.5-\xi} F_j(x) dx, \;\;\; {\rm for} \; v\leq 0, \\
\end{array}
 \right.\label{eq:4}\\
\xi&=&\frac{1}{\Delta x}\int_{t_{n}}^{t_{n+1}} v dt.\label{eq:6}
\end{eqnarray}
$\xi$ represents the travel distance normalized by the grid width from the upstream position at time $t=t_{n}$ to the cell boundary located at $x_j+\sgn(v)\Delta x/2$ at time $t=t_{n+1}$, integrated along the characteristic curve. $F_j(x)$ is the reconstruction function centered at $x_j$, which affects both the accuracy and stability of the solution.
\reviser{Although Equation (\ref{eq:4}) is presented for small CFL conditions $(|\xi| \le 1)$ for simplicity, the scheme can be extended to large CFL conditions by integrating the profile from the upwind position, which may lie outside the $j$-th cell, to the boundary of the $j$-th cell}.

\cite{2001JCoPh.172..166F} developed the PFC scheme that utilizes a piecewise quadratic polynomial for the reconstruction function,
\begin{eqnarray}
 F_j(x) &=& f_j^n - \frac{a_2}{12} + a_1 x + a_2 x^2 \nonumber \\
&=& \left[ \alpha_+f_j^n + \frac{a_1+a_2}{6}\left(3x^2+3x-\frac{1}{4}\right)\right] \nonumber \\
&& +  \left[ \alpha_- f_j^n - \frac{a_1-a_2}{6}\left(3x^2-3x-\frac{1}{4}\right)\right] \nonumber \\
&& + (1-\alpha_+-\alpha_-) f_j^n,\label{eq:5}
\end{eqnarray}
where the function meets the constraint $\int_{-0.5}^{0.5}F_j(x) dx=f_j^n$.
Formal second-order accuracy is attained by using optimal coefficients,
\begin{eqnarray}
 a_1 = \frac{f_{j+1}^n-f_{j-1}^n}{2}, \;\;\; a_2 = \frac{f_{j+1}^n-2f_j^n+f_{j-1}^n}{2}.\label{eq:21}
\end{eqnarray}
To ensure stability against numerical oscillations caused by the linear high-order reconstruction, the function is bounded by minimum and maximum values, denoted as $f_{j,\min} (\leq \min(f_{j-1}^n,f_{j}^n,f_{j+1}^n))$ and $f_{j,\max} (\geq \max(f_{j-1}^n,f_{j}^n,f_{j+1}^n))$, for $\forall x \in [-0.5,0.5]$.
This condition is satisfied provided that
\begin{eqnarray}
 \alpha_{\pm}\geq 0, \; 1-\alpha_+-\alpha_- \geq 0,\label{eq:8}\\
 \alpha_+ f_{j,\min} \leq \left[ \alpha_+f_j^n + \frac{a_1+a_2}{6}\left(3x^2+3x-\frac{1}{4}\right)\right] \leq \alpha_+ f_{j,\max},\label{eq:9}\\
 \alpha_- f_{j,\min} \leq \left[ \alpha_-f_j^n - \frac{a_1-a_2}{6}\left(3x^2-3x-\frac{1}{4}\right)\right] \leq \alpha_- f_{j,\max}.\label{eq:10}
\end{eqnarray}
Consequently, we derive the following inequalities for the coefficients $a_1$ and $a_2$,
\begin{eqnarray}
&& \alpha_{\pm} f_{j,{\rm min{\pm}}} \leq a_1{\pm}a_2 \leq \alpha_{\pm} f_{j,{\rm max{\pm}}},\label{eq:11}\\
&& f_{j,{\rm min+}}=3\max\left[2\left(f_j^n-f_{j,\max}\right),\left(f_{j,\min}-f_j^n\right)\right],\label{eq:27}\\
&& f_{j,{\rm max+}}=3\min\left[2\left(f_j^n-f_{j,\min}\right),\left(f_{j,\max}-f_j^n\right)\right],\label{eq:28}\\
&& f_{j,{\rm min-}}=3\max\left[2\left(f_{j,\min}-f_j^n\right),\left(f_j^n-f_{j,\max}\right)\right]=-f_{j,{\rm max+}},\label{eq:29}\\
&& f_{j,{\rm max-}}=3\min\left[2\left(f_{j,\max}-f_j^n\right),\left(f_j^n-f_{j,\min}\right)\right]=-f_{j,{\rm min+}},\label{eq:30}
\end{eqnarray}
\reviser{and the coefficients are corrected as follows,
\begin{eqnarray}
a_1 + a_2 \leftarrow {\rm median}(a_1 + a_2,\alpha_{+} f_{j,{\rm min{+}}},\alpha_{+} f_{j,{\rm max{+}}})={\rm median}(f_{j+1}^n-f_j^n,\alpha_{+} f_{j,{\rm min{+}}},\alpha_{+} f_{j,{\rm max{+}}}),\label{eq:66} \\
a_1 - a_2 \leftarrow {\rm median}(a_1 - a_2,\alpha_{-} f_{j,{\rm min{-}}},\alpha_{-} f_{j,{\rm max{-}}})={\rm median}(f_{j}^n-f_{j-1}^n,\alpha_{-} f_{j,{\rm min{-}}},\alpha_{-} f_{j,{\rm max{-}}}).\label{eq:67}
\end{eqnarray}
These equations reduce to the slope corrector introduced in the original study by setting $\alpha_+ = \alpha_- = 1/3$.}

While \cite{2001JCoPh.172..166F} employed the parameters $f_{j,\min}=0$, and $f_{j,\max}=\max(f_0^n,f_1^n,\dots)$, \cite{2008EP&S...60..773U} improved this method by implementing positivity-preserving and non-oscillatory bounds, which can be extended as follow,
\begin{eqnarray}
 f_{j,\max}&=&\max\left(f_{j-1/2,\max},f_{j+1/2,\max}\right),\label{eq:13}\\
f_{j,\min}&=&\max\left[0,\min\left(f_{j-1/2,\min},f_{j+1/2,\min}\right)\right],\label{eq:14}\\
f_{j-1/2,\max}&=&\max\left[\max\left(f_{j-1}^n,f_{j}^n\right),\min\left(f_{j-1/2,L},f_{j-1/2,R}\right)\right],\label{eq:15}\\
f_{j-1/2,\min}&=&\min\left[\min\left(f_{j-1}^n,f_{j}^n\right),\max\left(f_{j-1/2,L},f_{j-1/2,R}\right)\right],\label{eq:16}\\
f_{j-1/2,L}&=&f_{j-1}^n+r\left(f_{j-1}^n-f_{j-2}^n\right)+\left(1-r\right)\left(f_{j}^n-f_{j-1}^n\right),\label{eq:47}\\
f_{j-1/2,R}&=&f_{j}^n+r\left(f_{j}^n-f_{j+1}^n\right)+\left(1-r\right)\left(f_{j-1}^n-f_{j}^n\right),\label{eq:48}
\end{eqnarray}
where we employ $r=2/3$ derived from a quadratic polynomial, instead of $r=1$ used in the original study.

The resolution can be improved by setting large possible values to $\alpha_{\pm}$ that control the permissible range of the slope through Equations (\ref{eq:11}) to (\ref{eq:30}).
Selecting the original values of $\alpha_+ = \alpha_- = 1/3$ is a modest choice, and they can be increased to $\alpha_+ = \alpha_- = 1/2$.
Furthermore, defining $\alpha_{\pm}$ as a function of $a_1 \pm a_2$ is a reasonable approach to provide greater room for larger $|a_1 \pm a_2|$.
Consequently, we optimize $\alpha_{\pm}$ as
\begin{eqnarray}
\alpha_{\pm} &=& \frac{\beta_{\pm}}{\beta_+ + \beta_-},\label{eq:18}\\
\beta_{\pm} &=& \min\left(\beta_{\pm}^*+\epsilon,1\right),\label{eq:19}\\
\beta_{\pm}^*&=&\left\{
\begin{array}{l}
{(a_1{\pm}a_2)}/{(f_{j,{\rm max{\pm}}}+\epsilon)}, \;\;\; {\rm for} \; a_1{\pm}a_2>0,\\
{(a_1{\pm}a_2)}/{(f_{j,{\rm min{\pm}}}-\epsilon)}, \;\;\; {\rm for} \; a_1{\pm}a_2 \leq 0,
\end{array}
\right.\label{eq:31}
 \end{eqnarray}
where $\epsilon = 10^{-7}$ is a small positive constant to prevent division by zero.

Once the reconstruction function is determined, the outgoing flux in Equation (\ref{eq:4}) is calculated using Simpson's rule, e.g.,
\begin{eqnarray}
 \Phi_j = \int_{0.5-\xi}^{0.5} F_j(x) dx =  \frac{F_j(0.5)+4F_j(0.5(1-\xi))+F_j(0.5-\xi)}{6}\xi.\label{eq:20}
\end{eqnarray}

\section{Weighted positive and flux conservative (WPFC) method}\label{sec:weight-posit-flux}
The third-order PFC scheme is hard to simultaneously preserve extrema and non-oscillatory property, as a three-point stencil cannot make a distinction between an extremum and a discontinuity \citep{1997JCoPh.136...83S}. 
We extend the third-order PFC scheme to attain fifth-order accuracy, following an approach similar to the conservative semi-Lagrangian WENO scheme \citep{qiu2010conservative}.
The $j$-th stencil $S_j \in [j-2,j+2]$ is composed of three substencils: $S_{j,0} \in [j-2,j],S_{j,1} \in [j-1,j+1]$, and $S_{j,2} \in [j,j+2]$.
The piecewise quadratic polynomials are reconstructed in each substencil,
\begin{eqnarray}
 F_{j,k}(x) &=& f_j^n - \frac{a_{2,k}}{12} + a_{1,k} x + a_{2,k} x^2, \;\;\; k=0,1,2.\label{eq:22}
\end{eqnarray} 
Optimal coefficients for these polynomials are given by
\begin{eqnarray}
 a_{1,0}^{\rm opt} = \frac{3f_{j}^n-4f_{j-1}^n+f_{j-2}^n}{2}, \;\;\; a_{2,0}^{\rm opt} = \frac{f_{j}^n-2f_{j-1}^n+f_{j-2}^n}{2},\label{eq:23}\\
 a_{1,1}^{\rm opt} = \frac{f_{j+1}^n-f_{j-1}^n}{2}, \;\;\; a_{2,1}^{\rm opt} = \frac{f_{j+1}^n-2f_{j}^n+f_{j-1}^n}{2},\label{eq:24}\\
 a_{1,2}^{\rm opt} = \frac{-f_{j+2}^n+4f_{j+1}^n-3f_{j}^n}{2}, \;\;\; a_{2,2}^{\rm opt} = \frac{f_{j+2}^n-2f_{j+1}^n+f_{j}^n}{2}.\label{eq:25}
\end{eqnarray}
These coefficients are corrected to ensure positivity-preserving and non-oscillatory properties,
\begin{eqnarray}
 a_{1,k}{\pm}a_{2,k} = \median\left( a_{1,k}^{\rm opt} {\pm} a_{2,k}^{\rm opt}, \alpha_{\pm} f_{j,{\rm min{\pm}}}, \alpha_{\pm} f_{j,{\rm max{\pm}}}\right),\label{eq:26}
\end{eqnarray}
where \reviser{$f_{j,{\rm min{\pm}}},f_{j,{\rm max{\pm}}}$ are given by Equations (\ref{eq:27}) to (\ref{eq:30}) and (\ref{eq:13}) to (\ref{eq:48}),} and $\alpha_{\pm}$ is optimized by using $a_{1,k}^{\rm opt} {\pm} a_{2,k}^{\rm opt}$ in Equations (\ref{eq:18}) to (\ref{eq:31}).
The reconstruction function in $S_{j}$ is obtained from a convex combination of the polynomials in $S_{j,k}$,
\begin{eqnarray}
 F_j(x) = \sum_{k=0}^2w_k F_{j,k}(x).\label{eq:12}
\end{eqnarray}
Formal fifth-order accuracy is attained when the weights $w_k$ are set to their optimal weights $d_k$ for conservative semi-Lagrangian schemes,
\begin{eqnarray}
d_{1-\sgn(v)} = \frac{2+3|\xi|+\xi^2}{20},d_1=\frac{6+|\xi|-\xi^2}{10},d_{1+\sgn(v)}=\frac{6-5|\xi|+\xi^2}{20},\label{eq:32}
\end{eqnarray}
where $\xi$ is given in Equation (\ref{eq:6}). 
\cite{qiu2010conservative} designed the weights $w_k$ based on the methodology of the WENO schemes; assigning high weights to substencils with smooth profiles, low weights to those with discontinuous profiles, and weights close to the optimal values when the profile is sufficiently smooth within a stencil.
However, this approach is not suitable for the present scheme, because the profiles in the substencils are already positive and non-oscillatory.
Therefore, it is unnecessary to assign higher weights to smoother substencils to ensure numerical stability.
Instead, we can design the weights to outperform the underlying fifth-order scheme.

In the simulation of the advection equation, upwind schemes inherently introduce a diffusion term to stabilize the simulation such as
\begin{eqnarray}
\pdif{f}{t} + v \pdif{f}{x} = \pdif{}{x}\left(\nu_{\rm num} \pdif{f}{x}\right),\label{eq:33} 
\end{eqnarray}
where $\nu_{\rm num}$ represents the scheme-dependent numerical diffusion coefficient and $v$ is constant for simplicity.
By multiplying $f$ to this equation and then integrating it over space, diffusion consistently reduces the $L^2$ norm of the profile,
\begin{eqnarray}
\pdif{}{t} \int \frac{f^2}{2} dx = - \int \nu_{\rm num} \left( \pdif{f}{x}\right)^2 dx.\label{eq:34}
\end{eqnarray}
Based on this consideration, we anticipate that numerical diffusion can be improved by designing the scheme to suppress the decrease in the $L^2$ norm.
To achieve this, we formulate the nonlinear weights $w_k$ that assign higher priority to substencils with larger $L^2$ norm.
For quadratic polynomials (Eq. (\ref{eq:22})), the piecewise $L^2$ norm is calculated as
\begin{eqnarray}
 L_{j,k}^2 = \int_{-0.5}^{0.5} F_{j,k}^2(x) dx= (f_j^n)^2+\frac{a_{1,k}^2}{12} + \frac{a_{2,k}^2}{180},\label{eq:35}
\end{eqnarray}
which indicates that steeper substencils have larger $L^2$ norm.
In that regard, our approach is similar to that of the improved WENO-Z schemes \citep{acker2016improved,luo2021improved}, in which weights for less smooth substencils are increased to enhance resolution from the underlying WENO-Z scheme \citep{2008JCoPh.227.3191B}.

We define the nonlinear weights as
\begin{eqnarray}
w_k = \frac{\gamma_k}{\sum_{k=0}^{2} \gamma_k}, \;
 \gamma_k = d_k \left[C + \left(\frac{\Delta L_{j,k}^2 +\epsilon}{\Delta L_{j}^2 +\epsilon}\right)^p \right],\label{eq:36}\\
\Delta L_{j,k}^2 = L_{j,k}^2 - (f_j^n)^2 = \frac{a_{1,k}^2}{12} + \frac{a_{2,k}^2}{180},\label{eq:40}\\
\Delta L_{j}^2 = L_{j}^2 - (f_j^n)^2 = \frac{a_{1}^2}{12} + \frac{a_{2}^2}{180} + \frac{a_3^2}{448} + \frac{a_4^2}{3600} + \frac{a_1 a_3}{40} + \frac{a_2 a_4}{420},\label{eq:41}  
\end{eqnarray}
where $\epsilon=10^{-7}$ and $L_{j}^2$ in Equation (\ref{eq:41}) is the piecewise $L^2$ norm for the fourth-degree polynomial in $S_j$,
\begin{eqnarray}
 F_{j,{\rm 4th}}(x) = f_j^n-\frac{a_2}{12}-\frac{a_4}{80}+a_1 x+a_2 x^2 + a_3 x^3 + a_4 x^4,\label{eq:42}\\
 a_1 = \frac{-5f_{j+2}^n+34f_{j+1}^n-34f_{j-1}^n+5f_{j-2}^n}{48} = 
\frac{5a_{1,0}^{\rm opt}+14a_{1,1}^{\rm opt}+5a_{1,2}^{\rm opt}}{24},\label{eq:43}\\
 a_2 = \frac{-f_{j+2}^n+12f_{j+1}^n-22f_j^n+12f_{j-1}^n-f_{j-2}^n}{16} =
\frac{-a_{2,0}^{\rm opt}+10a_{2,1}^{\rm opt}-a_{2,2}^{\rm opt}}{8},\label{eq:44}\\
 a_3 = \frac{f_{j+2}^n-2f_{j+1}^n+2f_{j-1}^n-f_{j-2}^n}{12} = 
\frac{-a_{1,0}^{\rm opt}+2a_{1,1}^{\rm opt}-a_{1,2}^{\rm opt}}{6},\label{eq:45}\\
 a_4 = \frac{f_{j+2}^n-4f_{j+1}^n+6f_j^n-4f_{j-1}^n+f_{j-2}^n}{24} = 
\frac{a_{2,0}^{\rm opt}-2a_{2,1}^{\rm opt}+a_{2,2}^{\rm opt}}{12}.\label{eq:46}
\end{eqnarray}
The constant $C$ and the exponent $p$ in Equation (\ref{eq:36}) are positive tunable parameters that should satisfy the following criteria:
 \begin{enumerate}
  \item Set small and large possible values as long as ensuring that the numerical growth rate remains negative for all wavenumbers.
\item Computationally efficient.
 \end{enumerate}
Based on the approximate dispersion relation analysis conducted in the next subsection, we find that $C=p=1/2$ provides satisfactory results.
Using these weights, the reconstruction function and the outgoing flux are calculated through Equations (\ref{eq:12}) and (\ref{eq:20}), thereby completing the time integration.

\subsection{Approximate dispersion relation}\label{sec:appr-disp-relat}
\begin{figure}[htbp]
\centering
\includegraphics[clip,angle=0,scale=0.6]{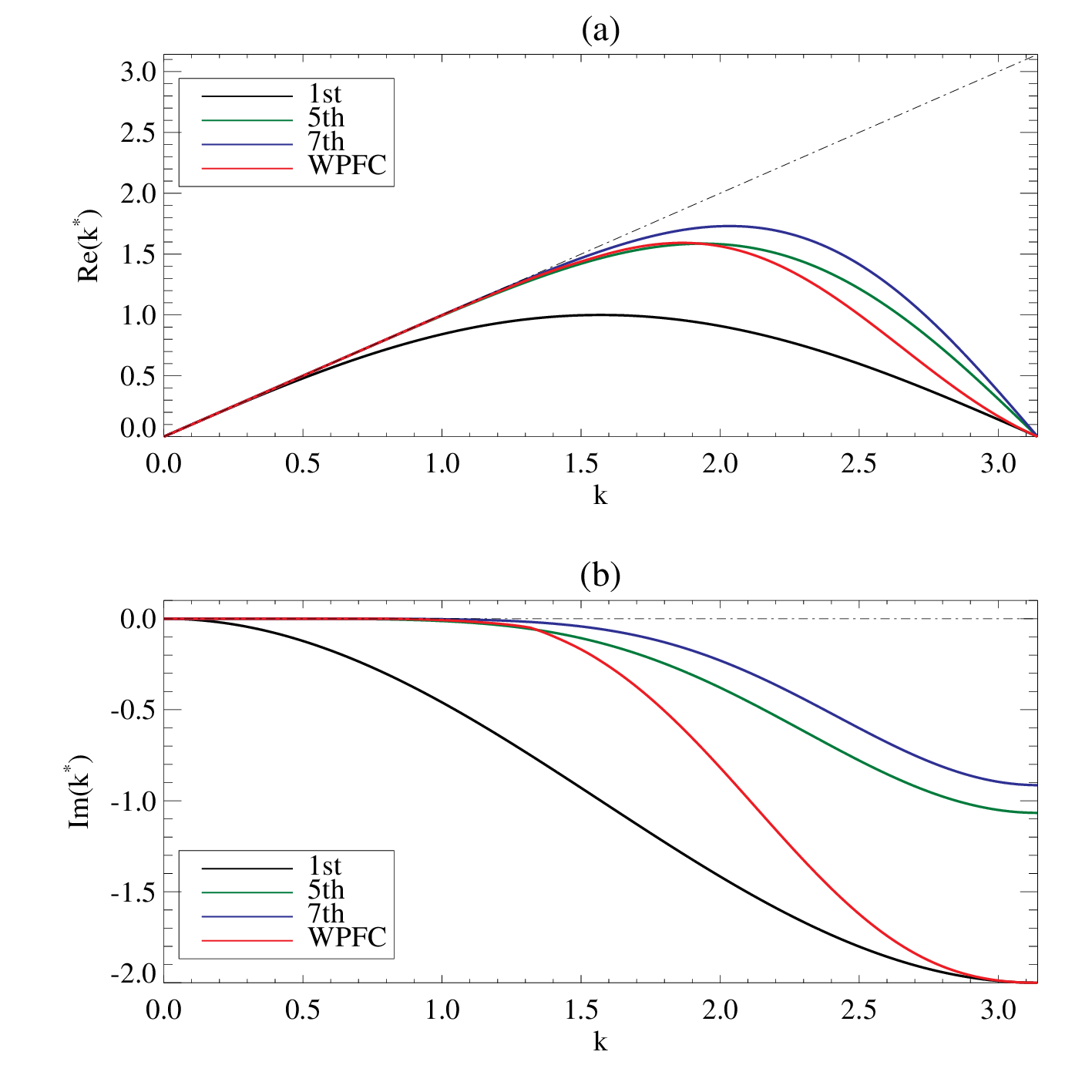}
\caption{Approximate dispersion relation. (a) Real and (b) imaginary parts of the modified wavenumber in the present scheme are shown in red lines. For comparison, results from linear first-, fifth-, and seventh-order schemes are shown in black, green, and blue lines. Dash-dotted lines indicate the exact solution.}
\label{fig:adr} 
\end{figure}

Following \cite{1992JCoPh.103...16L} and \cite{pirozzoli2006spectral}, we investigate the dispersion relation of the present scheme to assess its spectral properties.
By using the periodic function $f_j=\exp(i k x_j)$ and $\xi=0$, the modified wavenumber is calculated as
\begin{eqnarray}
 k_j^* = \frac{F_{j}(0.5)-F_{j-1}(0.5)}{i f_j \Delta x},\label{eq:37}
\end{eqnarray}
and then is averaged over space.
Figure \ref{fig:adr} shows the real and imaginary parts of the modified wavenumber as a function of the true wavenumber.
For wavenumbers below $k \sim 1.34$, both the numerical dispersion (real part) and dissipation (imaginary part) of the present scheme are better than those of the underlying linear fifth-order scheme.
Above $k \sim 1.34$, the numerical dissipation rate of the present scheme increases due to the inclusion of the slope corrector (Eq. (\ref{eq:26})). 
As $k$ approaches $\pi$, the dissipation rate becomes nearly equivalent to that of the first-order scheme.
Without the correction, the scheme produces oscillatory solutions similar to those from linear high-order schemes.
The imaginary part remains negative for all wavenumber range.
Its maximum value becomes positive \reviser{if $C < 1/2$ or $p \geq 1$ in Equation (\ref{eq:36})}, which is an unfavorable situation for numerical stability.

\begin{figure}[htbp]
\centering
\includegraphics[clip,angle=0,scale=0.6]{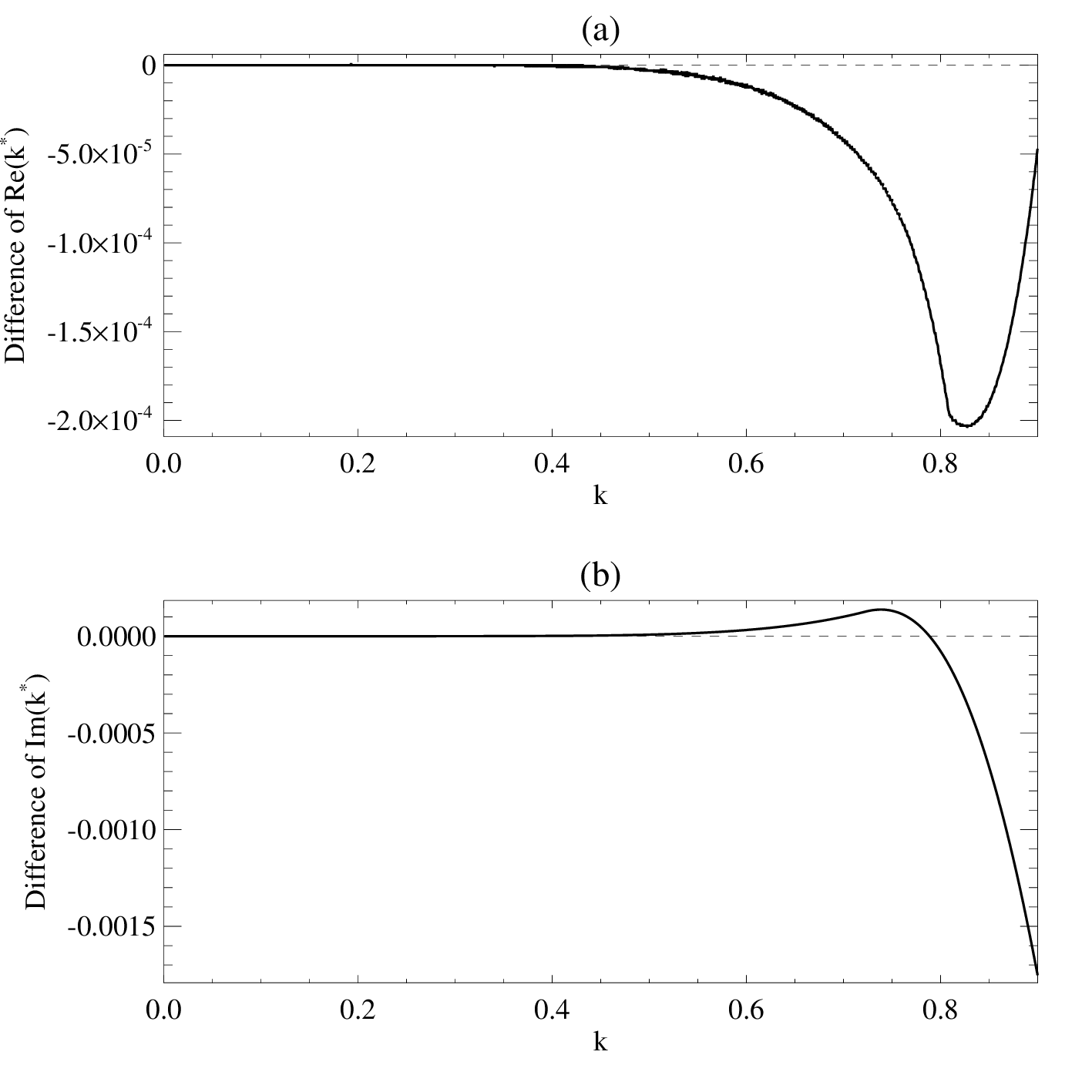}
\caption{Difference in (a) the dispersion error and (b) the dissipation error between the present scheme and the linear seventh-order scheme. The spectral properties of the present scheme are superior when the difference is negative in (a) and positive in (b), respectively.}
\label{fig:dadr} 
\end{figure}

It is worth noting that the spectral properties of the present scheme even surpass those of the seventh-order scheme in certain wavenumber ranges.
Figure \ref{fig:dadr} shows the difference in the dispersion and dissipation errors between the present scheme and the linear seventh-order scheme,
$ |Re(k_{\rm present}^*)-k| - |Re(k_{\rm 7th}^*)-k|$, and $Im(k_{\rm present}^*)-Im(k_{\rm 7th}^*).$
For wavenumbers in the range $0.062 < k < 0.79$, the present scheme exhibits lower dispersion error and dissipation rate compared to the seventh-order scheme.
These findings indicate that well-designed nonlinear weighting can enhance both the accuracy and stability of the scheme, potentially outperforming not only the underlying scheme but also higher-order schemes.

\section{Numerical experiments}\label{sec:numer-exper}
We conduct numerical experiments of linear advection and Vlasov equations to assess the capability of the present scheme.

\subsection{Linear advection}\label{sec:linear-advection}
\begin{figure}[htbp]
\centering
\includegraphics[clip,angle=0,scale=0.6]{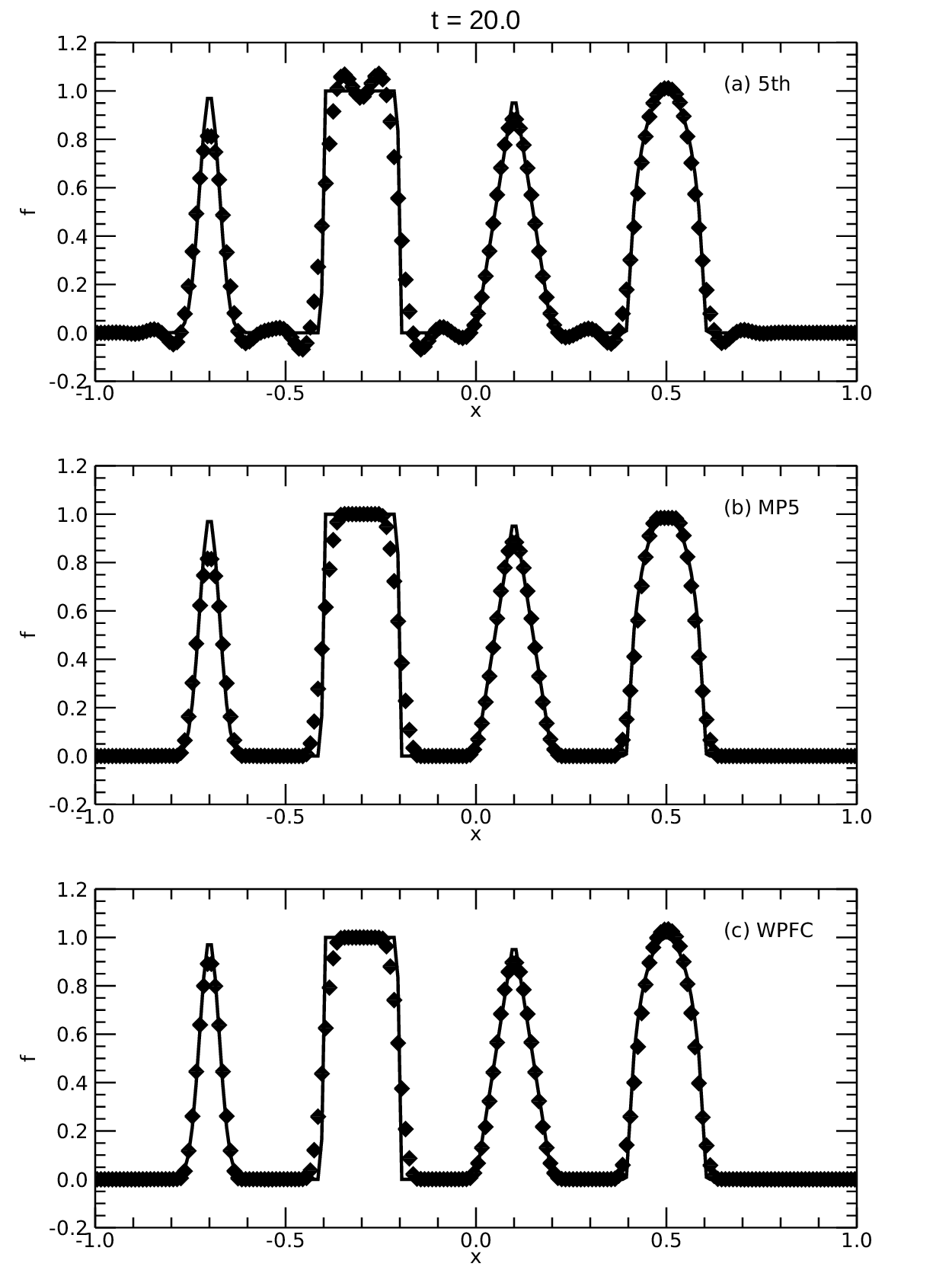}
\caption{Linear advection over 10 periods obtained from fifth-order conservative semi-Lagrangian schemes: (a) the linear scheme, (b) the monotonicity-preserving scheme, and (c) the present scheme. Solid lines indicate the exact solution.}
\label{fig:chk5} 
\end{figure}

We solve Equation (\ref{eq:1}) using $v=1$ and a CFL number of 0.4.
The spatial domain spans $[-1,1]$ with the periodic boundary condition.
For comparison, we employ three different fifth-order conservative semi-Lagrangian schemes: the linear scheme, the monotonicity-preserving (MP5) scheme \citep{2017ApJ...849...76T}, and the present scheme.
Figure \ref{fig:chk5} shows the numerical solutions at $t=20$ for a combination of Gaussian, rectangular, triangle, and half ellipse waves resolved by 200 grid points, employed in previous studies \citep[e.g.,][]{1996JCoPh.126..202J,1997JCoPh.136...83S,2019ApJS..242...14M}.
The MP5 and present schemes completely eliminate the numerical oscillations observed in the linear scheme, and preserve positivity.
Differences between the MP5 and present schemes can be found at the peaks of the Gaussian wave (leftmost) and the half ellipse wave (rightmost).
Notably, the Gaussian peak in the present scheme is better preserved, indicating superior resolution compared to the linear and MP5 schemes.

Table \ref{tab:acc1} examines the numerical errors in the $L^1$ and $L^{\infty}$ norms, and the order of accuracy measured at $t=4.0$ for the Gaussian profile, $f(x)=\exp(-x^2/2s^2)$ where $s=1/16$. 
All three schemes show convergent results, and the order of accuracy approaches their formal fifth order as the number of grid points increases to 256.
In the present scheme, both $L^1$ and $L^{\infty}$ errors are approximately 2 times smaller than those in the linear and MP5 schemes.
Table \ref{tab:acc2} shows the results for the sinusoidal profile, $f(x)=[3+\sin(4 \pi x)]/4$.
For this profile, the MP5 scheme yields results identical to the linear scheme.
The present scheme shows the numerical errors that are an order of magnitude smaller than those in the linear scheme, and the order of accuracy measured by the $L^1$ error exceeds the formal fifth order for grid points ranging from 64 to 256.
These results confirm the superior performance of the present scheme compared to the underlying fifth-order scheme.

\reviser{The stability of the solution is sensitive to the choice of the parameters $C$ and $p$ in Equation (\ref{eq:36}).
The nonlinear weight $\gamma_k$ increases the priority of the steep substencil as $(\Delta L_{j,k}^2)^p/C$ increases.
Too small $C$ or too large $p$ cause a leading dispersion error and a growing dissipation error in the range $k \le 1.5$ in the approximate dispersion relation (not shown).
Consequently, the solution of the advection equation is significantly degraded, exhibiting overshoot and staircasing, although the positivity is still preserved.}

To extend the one-dimensional scheme to higher dimensions, we use time splitting methods \citep[e.g.,][]{2003CoPhC.150..247F}.
While the splitting method simplifies the application of the one-dimensional scheme by solving each dimension sequentially, it can degrade the overall accuracy of the solution due to low-order splitting in time.
To confirm that the application of high-order one-dimensional schemes is advantageous in multidimensional problems, we conduct the three-dimensional solid body rotation problem, which appears in electromagnetic Vlasov-Maxwell simulations,
\begin{eqnarray}
\pdif{f}{t}+\left(\vect{x} \times \vect{\Omega}\right) \cdot \nabla f=0,\label{eq:59}
\end{eqnarray}
where $\vect{\Omega} = (2 \pi/\mart{6},2 \pi/\mart{3},2 \pi/\mart{2})$.
The spatial domain spans $[-0.5,0.5]^3$ with the open boundary condition.
We employ an efficient splitting method proposed by \cite{2006CoPhC.175...86S}, \reviser{which builds upon the method of characteristics in conjunction with the Boris scheme \citep{PIC}.
The time step is set to be $\Delta t = 1/(10N_x)$, yielding a CFL number of ${\rm max}[(\vect{x} \times \vect{\Omega})_x/\Delta x,(\vect{x} \times \vect{\Omega})_y/\Delta y,(\vect{x} \times \vect{\Omega})_z/\Delta z] \Delta t=0.39$, where $N_x=N_y=N_z$ is the number of grid points in each dimension.}

Table \ref{tab:acc3} presents the numerical errors in the $L^1$ and $L^{\infty}$ norms measured after 10 rotations for the Gaussian profile, $f(x,y,z)=\exp(-x^2/2s_x^2-y^2/2s_y^2-z^2/2s_z^2)$ where $(s_x,s_y,s_z)=(0.06,0.08,0.1)$, solved by the linear fifth-order, MP5, and present schemes.
All schemes show convergence, and the order of accuracy of the schemes approaches their formal fifth order with increasing the number of grid points to $128^3$, indicating that the numerical error is primarily dominated by spatial discretization rather than temporal discretization.
This result demonstrates the benefit of utilizing high-order schemes.
Similar to the one-dimensional problems, the present scheme achieves smaller errors than the linear and MP5 schemes for grid points ranging from $32^3$ to $128^3$.
Thanks to the splitting method to apply the one-dimensional scheme dimension by dimension, positivity is preserved in multidimensional problems.

\reviser{Since the Boris scheme is second-order in time, temporal discretization will eventually dominate the numerical error as the spatial grid is refined while keeping $\Delta t/\Delta x$ constant.
The temporal discretization error can be minimized by using smaller time steps or an exact splitting method \citep[e.g.,][]{2025JCoPh.52913858L}.
To confirm this effect, we conduct the same problem using $\Delta t \propto (\Delta x)^{2.5}$.
The results are presented in Table \ref{tab:acc4}, demonstrating improved solution accuracy only on the finest grid in comparison with Table \ref{tab:acc3}.
Therefore, we consider the present approach to be sufficient for this study.
}

\subsection{Vlasov-\Ampere simulation}\label{sec:vlas-poiss-simul}
We solve one-dimensional electrostatic Vlasov-\Ampere equations as a demonstration of kinetic plasma simulations,
\begin{eqnarray}
&& \pdif{f}{t}+v\pdif{f}{x}+\frac{qE}{m}\pdif{f}{v}=0,\label{eq:56}\\
&& \pdif{E}{t}=-4 \pi q \int vfdv,\label{eq:57}
\end{eqnarray}
where $f(x,v,t)$ represents the phase space distribution function of electrons (ambient ions are assumed to be stationary), $q=-e$ and $m$ are elementary charge and mass of electrons, and $E(x,t)$ is the electric field, respectively.
\revrr{The electric field should satisfy Gauss's law as a constraint,}
\begin{eqnarray}
\pdif{E}{x} = 4 \pi e \left(1-\int fdv\right),\label{eq:65} 
\end{eqnarray}
\revrr{where the density is normalized to unity.}
The velocity space spanning $[-7v_t,7v_t]$ is discretized by $N_v = 64$, $128$ (fiducial), or $256$ grid points, where $v_t=1$ is the thermal velocity.
Open and periodic boundary conditions are imposed in the $v$- and $x$-direction.
The time step is set to $0.01 \omega_{pe}^{-1}$, where $\omega_{pe}=1$ is the electron plasma frequency.
The system of equations is updated using the standard operator splitting technique as follows: 
(1) Using $E$ at $n$ steps, solve $\partial f/\partial t + (qE/m)\partial f/\partial v = 0$ to advance $f$ in velocity space from $n-1/2$ to $n+1/2$ steps. 
(2) Solve $\partial f/\partial t + v\partial f/\partial x = 0$ to advance $f$ in physical space from $n$ to $n+1$ steps \revrr{using a conservative scheme. The numerical flux is defined as $vf_{j+1/2}$ where $j+1/2$ is the half index in physical space.} 
\revrr{(3) Integrate this numerical flux over velocity space to calculate the electric current in the right hand side of Equation (\ref{eq:57}).
Integrating the advection equation in Step (2) over velocity space, one find that this electric current satisfies the charge conservation law,}
\begin{eqnarray}
\pdif{}{t}q\int f_jdv = - \frac{q}{\Delta x}\left(\int vf_{j+1/2}dv - \int vf_{j-1/2}dv\right).\label{eq:68}
\end{eqnarray} 
(4) Solve Equation (\ref{eq:57}) to advance $E$ from $n$ to $n+1$ steps.
\revrr{By discretizing Equation (\ref{eq:57}) in physical space and combining it with Equation (\ref{eq:68}), Equation (\ref{eq:65}) is guaranteed to hold at the discrete level, provided that it is initially satisfied,}
\begin{eqnarray}
\pdif{}{t}\frac{E_{j+1/2}-E_{j-1/2}}{\Delta x}=\frac{-4 \pi q}{\Delta x} \left(\int vf_{j+1/2}dv - \int vf_{j-1/2}dv\right)=-4 \pi e \pdif{}{t}\int f_jdv.\label{eq:69}
\end{eqnarray}
\revrr{As a result, the scheme avoid a numerical cost to solve the Poisson equation for the electric potential \cite[e.g.,][]{2009JCoPh.228.4773S}.}
To investigate the effect of velocity space resolution on the quality of numerical solutions, the advection equation in velocity space (Step (1)) is updated by four different schemes: the third-order PFC scheme (in \S \ref{sec:posit-flux-cons}), fifth- and seventh-order MP schemes (MP5 and MP7, \cite{2017ApJ...849...76T}), and the present fifth-order weighted PFC scheme (WPFC, in \S \ref{sec:weight-posit-flux}).
The advection equation in physical space (Step (2)) is updated by the PFC scheme for simplicity, which employs the following compact positivity-preserving and non-oscillatory bounds,
\begin{eqnarray}
f_{j,\max}=\max\left[\max\left(f_{j-1}^n,f_{j+1}^n\right),
\min\left(2f_{j}^n-f_{j-1}^n,2f_{j}^n-f_{j+1}^n\right)\right],\label{eq:38}\\
f_{j,\min}=\max\left[0,\min\left\{\min\left(f_{j-1}^n,f_{j+1}^n\right),
\max\left(2f_{j}^n-f_{j-1}^n,2f_{j}^n-f_{j+1}^n\right)\right\}\right].\label{eq:39}
\end{eqnarray}

\begin{figure}[htbp]
\centering
\includegraphics[clip,angle=0,scale=0.25]{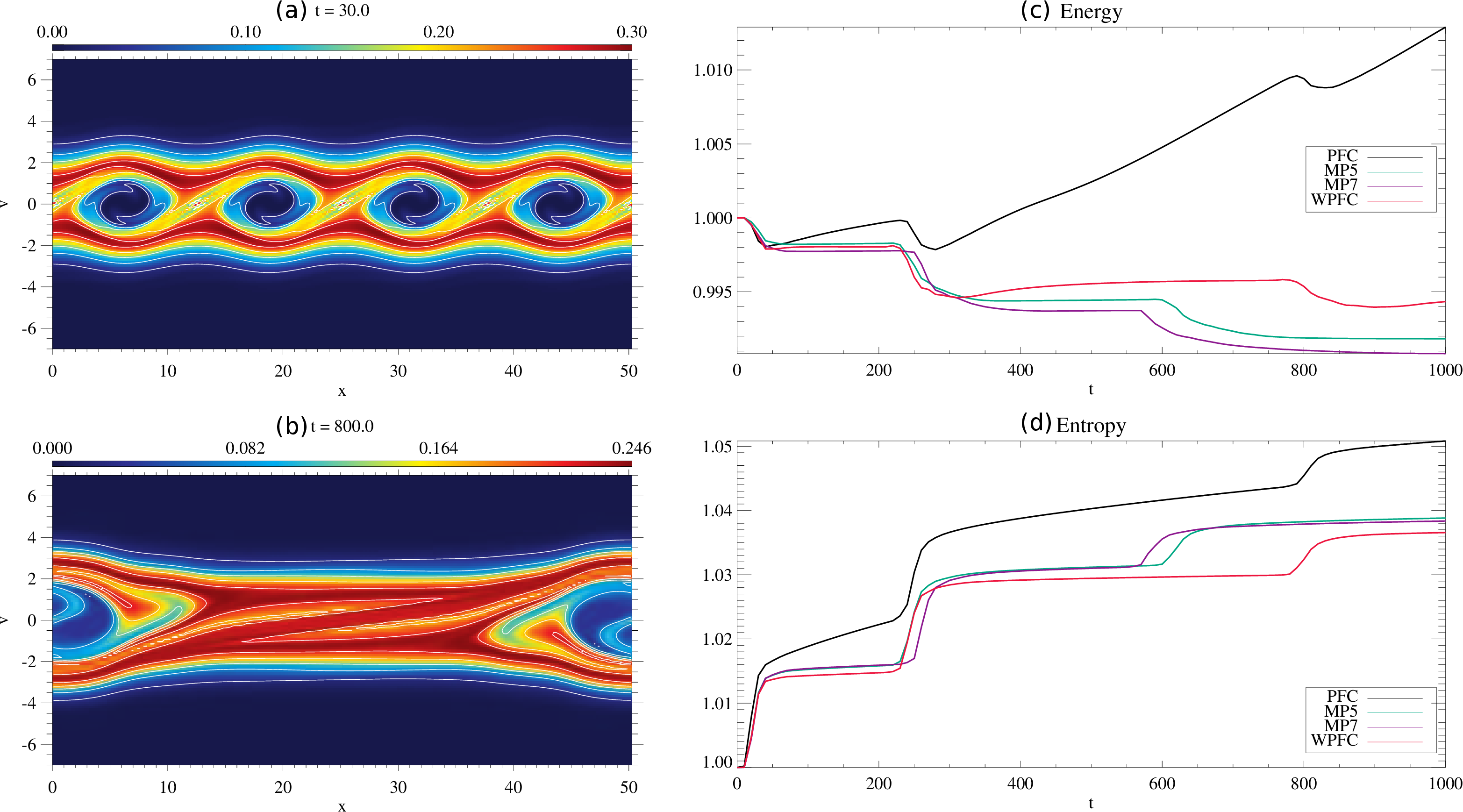}
\caption{Numerical experiment of the two stream instability. (a,b) Electron phase space distribution at $t=30 \omega_{\rm pe}^{-1}$ and $800 \omega_{\rm pe}^{-1}$. (c,d) Time profiles of the total energy and entropy normalized by their initial values, obtained from the third-order PFC scheme (black lines), the fifth- and seventh-order monotonicity-preserving schemes (green and purple lines), and the present scheme (red lines), respectively.}
\label{fig:twost} 
\end{figure}

The first problem is the two stream instability.
We set 256 grid points in physical space, and the initial condition of the phase space distribution as
\begin{eqnarray}
 f(v,x)=\frac{1 + 0.05 \cos \left(k x \right)}{\mart{2 \pi} v_t} \frac{v^2}{v_t^2} \exp\left( -\frac{v^2}{2 v_t^2} \right), x \in [0,16 \pi d_e],\label{eq:51}
\end{eqnarray}
 where $d_e=1$ is the Debye length and $k=0.5/d_e$ gives the mode number of the initial perturbation of $n=4$.
As the instability grows, electron hole structures emerge in phase space, as shown in Figure \ref{fig:twost}(a).
The solution agrees with previous studies \citep[e.g.,][]{2011JCoPh.230.6800M}, and with solutions obtained from other high-order conservative semi-Lagrangian schemes, including the WENO scheme \citep{qiu2010conservative} and the MP scheme.
Subsequently, the phase space holes merge, eventually forming the $n=1$ mode, as illustrated in Figure \ref{fig:twost}(b).

We assess the solution in terms of the conservation of total energy and entropy,
\begin{eqnarray}
 \left(\text{Energy}\right)  &=& \int\left(\int \frac{m v^2}{2}f dv + \frac{E^2}{8 \pi}\right) dx,\label{eq:53}\\
 \left(\text{Entropy}\right) &=& \int \!\!\! \int f (1-f) dv dx,\label{eq:54}
\end{eqnarray}
where Equation (\ref{eq:54}) is introduced by \cite{1988JSP....52..479T} as a generalized form for entropy.
Because the schemes conserve the $L^1$ norm, the increase in entropy corresponds to the decrease in the $L^2$ norm caused by numerical diffusion (Eq. (\ref{eq:34})).
Figure \ref{fig:twost}(c) and (d) show the time profiles of energy and entropy.
After the initial jump at $t=20\omega_{pe}^{-1}$ caused by the growth of the $n=4$ mode, both energy and entropy change suddenly as the phase space holes coalesce.
They are conserved within 1\% and 4\% in the MP and present schemes, respectively, whereas the PFC scheme shows a continuous increase over time.
All schemes demonstrate the improvement of conservation with increasing resolution, confirming the convergence of numerical solutions.
For the case of $N_v = 128$ presented here, the present scheme exhibits better conservation of energy compared to the MP schemes.
However, this trend does not always hold true for different $N_v$ values, complicating the assessment of quality based on energy conservation.
Meanwhile, the present scheme consistently exhibits better conservation of entropy compared to the other schemes for all cases.

The second problem is the bump-on-tail instability.
We set 256 grid points in physical space, and the initial condition of the phase space distribution composed of core and beam components as follows:
\begin{eqnarray}
 f(v,x)=\frac{1 + 0.05 \cos \left(k x\right)}{\mart{2 \pi}} \left[\frac{n_p}{v_t}\exp\left( -\frac{v^2}{2 v_t^2} \right) + \frac{n_b}{v_{t,b}}\exp\left( -\frac{(v-v_b)^2}{2 v_{t,b}^2} \right) \right],\nonumber \\ x \in [0,40 \pi d_e],\label{eq:58}
\end{eqnarray}
where $k=0.3/d_e$, $n_p=0.97$ and $n_b=0.03$ are core and beam densities, $v_{b}=4.5 v_{t}$ and $v_{t,b}=0.5 v_{t}$ are bulk and thermal velocities of the beam component, respectively.
The velocity distribution function in Equation (\ref{eq:58}) has a positive slope at a velocity corresponding to the phase velocity of the Langmuir wave, $v_{ph} = \mart{\omega_{pe}^2+3 k^2 d_e^2}/k = 3.756 v_t$, leading to the instability.
In this simulation, the background electric current $(qn_b v_b)$ is subtracted from Equation (\ref{eq:57}) to remove the homogeneous electric field component.
The result is equivalent to that obtained by solving the Poisson equation for the electric potential.

\begin{figure}[htbp]
\centering
\includegraphics[clip,angle=0,scale=0.25]{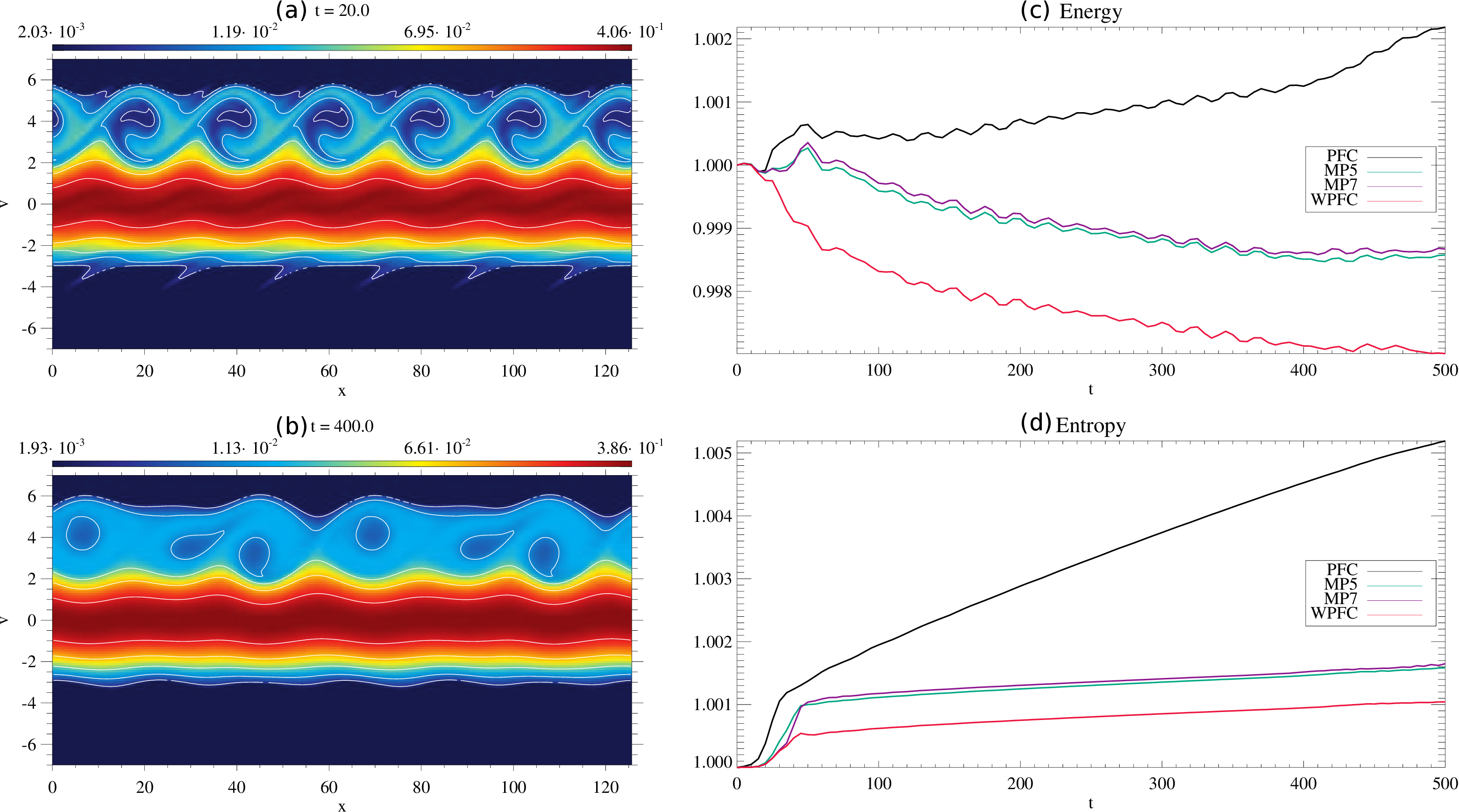}
\caption{Numerical experiment of the bump-on-tail instability. The figure follows the same format as Figure \ref{fig:twost}.}
\label{fig:buota} 
\end{figure}

Figure \ref{fig:buota}(a) shows the phase space distribution at a nonlinear stage, $t=20 \omega_{pe}^{-1}$.
The initial $n=6$ mode grows, forming phase space holes centered around $v=v_{ph}$.
These phase space holes remain stable until coalescence occurs at $t=400 \omega_{pe}^{-1}$ as shown in Figure \ref{fig:buota}(b) \citep[e.g.,][]{2017LPB....35..706S}.
This coalescence is resolved by the MP and present schemes for all cases with $N_v=64,128,258$, whereas the PFC scheme is difficult to resolve it at $N_v=64$.
Figure \ref{fig:buota}(c) and (d) show the time profiles of energy and entropy.
The MP and present schemes conserve entropy within 0.17\% and 0.1\%, respectively, whereas it increases rapidly over time in the PFC scheme.
For the case of $N_v = 128$ presented here, the MP schemes exhibit better conservation of energy compared to the present scheme.
Similar to the two stream instability problem, however, this trend does not always hold true for different $N_v$ values.
Meanwhile, the present scheme consistently exhibits better conservation of entropy compared to the other schemes for all cases.

\begin{figure}[htbp]
\centering
\includegraphics[clip,angle=0,scale=0.35]{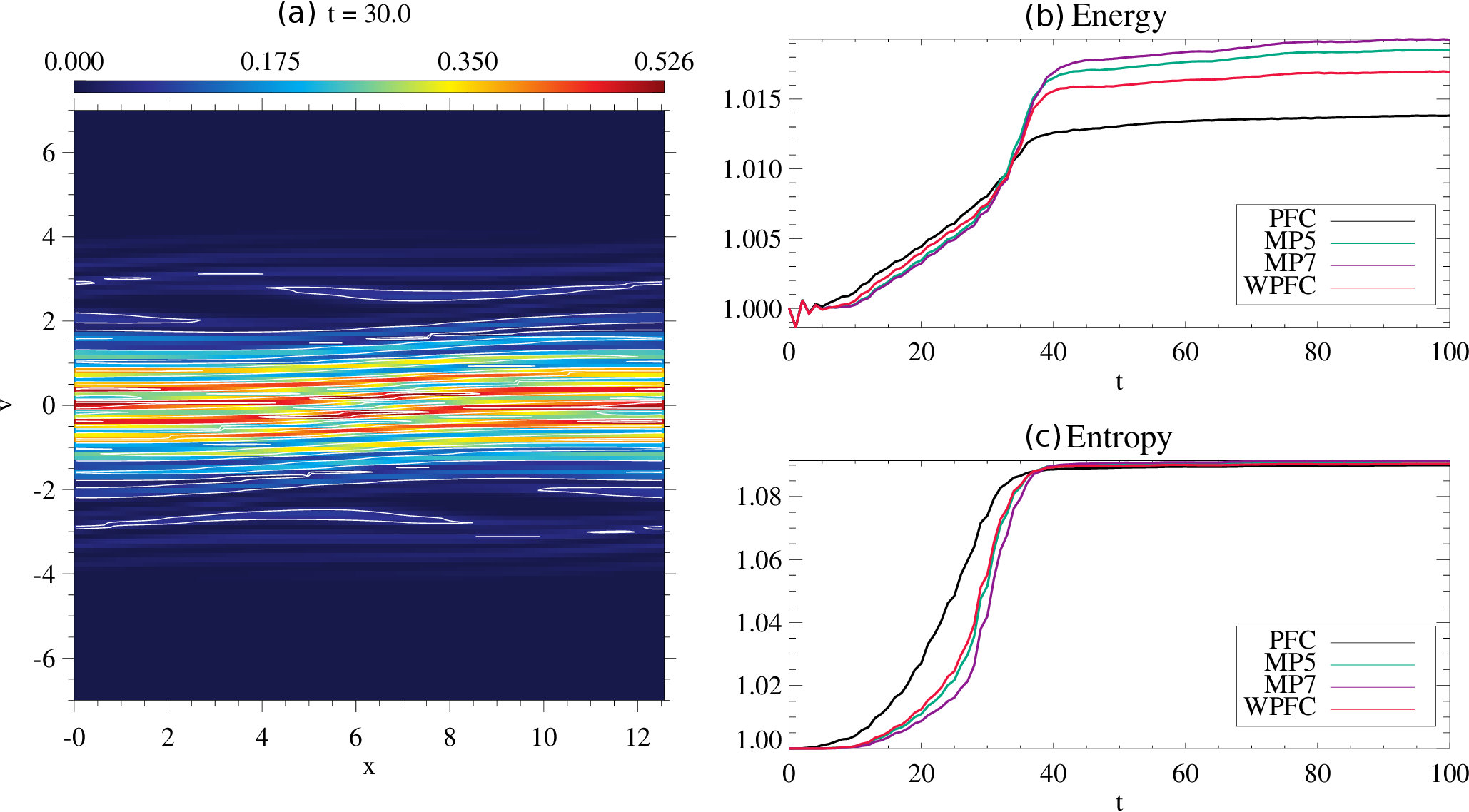}
\caption{Numerical experiment of the nonlinear Landau damping. The figure follows the same format as Figure \ref{fig:twost}. }
\label{fig:nllan} 
\end{figure}

The third problem is the nonlinear Landau damping.
We set 128 grid points in physical space, and the initial condition of the phase space distribution as
\begin{eqnarray}
 f(v,x)=\frac{1 + 0.5 \cos \left(k x\right)}{\mart{2 \pi} v_t} \exp\left( -\frac{v^2}{2 v_t^2} \right), x \in [0,4 \pi d_e],\label{eq:55}
\end{eqnarray}
where $k=0.5/d_e$.
Due to large inhomogeneity in physical space, particle streaming results in the formation of a pronounced filamentary structure in velocity space, as shown in Figure \ref{fig:nllan}(a).
The initial perturbation is damped by $t=15\omega_{pe}^{-1}$, followed by the growth of electric fields until $t=40\omega_{pe}^{-1}$, consistent with previous studies \citep[e.g.,][]{1976JCoPh..22..330C,1999CoPhC.120..122N}.
As the system evolves, the filamentary structure becomes progressively finer and eventually approaches the grid size.
The timing of dissipation of this structure is proportional to the number of grid points in velocity space.
Figure \ref{fig:nllan}(b) and (c) show the time evolution of energy and entropy.
Entropy begins to increase earlier in the PFC scheme than in the other schemes.
Eventually, the entropy increases by approximately 9\%, and is almost the same level among the schemes.
This result suggests that grid-scale dissipation in velocity space is a major cause of the numerical error of entropy. 
Increasing $N_v$ from 128 to 256 delays the onset of entropy increase, yet only marginally improves its saturation level by approximately 0.2\%.
Meanwhile, the conservation of energy varies across the schemes.
The PFC scheme exhibits the best energy conservation, followed by the present scheme.
This trend is observed at $N_v=64$ and $256$ as well, indicating that utilizing high-order schemes does not always improve conservation properties, at least for the nonlinear Landau damping problem, where grid-scale structure dominates the numerical error.

\subsection{Vlasov-Maxwell simulation}\label{sec:vlas-maxw-simul}
We solve two-dimensional electromagnetic Vlasov-Maxwell equations for ions and electrons,
\begin{eqnarray}
&& \pdif{f_s}{t}+\vect{v}\cdot \nabla_x f_s + \frac{q_s}{m_s} \left(\vect{E}+\frac{\vect{v} \times \vect{B}}{c}\right) \nabla_v f_s = 0,\;\;\; (s=i,e)\label{eq:60}\\
&& \pdif{\vect{B}}{t} = -c \nabla \times \vect{E},\label{eq:62}\\
&& \pdif{\vect{E}}{t} = c \nabla \times \vect{B} - 4 \pi \sum_{s=i,e} q_s \int \vect{v} f_s d\vect{v},\label{eq:61}
\end{eqnarray}
with constraints,
\begin{eqnarray}
\nabla \cdot \vect{B} &=& 0,\label{eq:64}\\
\nabla \cdot \vect{E} &=& 4 \pi \sum_{s=i,e} q_s \int f_s d\vect{v},\label{eq:63}
\end{eqnarray}
where $q_s$ and $m_s$ are elementary charge and mass of $s$ species ($i$ for ions and $e$ for electrons), $c$ is the speed of light, and $B(\vect{x},t)$ is the magnetic field, respectively.
Similar to the Vlasov-\Ampere simulations, the system of equations is advanced using the operator splitting technique, where the advection equation in velocity space is updated by the PFC, MP5, MP7, and present schemes to investigate the effect of velocity space resolution on the quality of numerical solutions.
The Maxwell equations (\ref{eq:62}) and (\ref{eq:61}) are updated by the standard finite-difference time-domain (FDTD) method, which guarantees the divergence-free condition of the magnetic field at the discrete level.
Gauss's law is also guaranteed by the same manner described in \S~\ref{sec:vlas-poiss-simul}.
The velocity space spanning $[-5v_{t,s},5v_{t,s}]^3$ is discretized by $N_v^3 = 32^3 $ grid points with imposing open boundary condition.

We test collisionless magnetic reconnection as demonstrated in \reviser{\cite{2006PhPl...13i2309S}, \cite{2009CoPhC.180..365U} and} \cite{2015CoPhC.187..137M}.
\reviser{Similar to the Geospace Environment Modeling magnetic reconnection challenge by \cite{2001JGR...106.3737B} and the study by \cite{2011PhPl...18l2108Z}, we adopt the Harris current sheet configuration with a stationary background plasma as the initial condition, $\vect{B}(y)=B_0 {\rm tanh}\left(y/\lambda\right) \vect{e}_x$ and $n(y)=n_0 {\rm cosh}^{-2}\left(y/\lambda\right)+n_{\rm bg}$,
where $n(y)$ is the density, $\lambda$ is the half thickness of the current sheet,} and the upstream magnetic pressure balances the plasma pressure within the  current sheet (see \cite{2015CoPhC.187..137M} for details). 
The ion-to-electron mass and temperature ratios are set to 25 and 1, respectively.
The ratios of the upstream electron gyrofrequency to the plasma frequency in the current sheet and the upstream density to the current sheet density are 0.1 and 0.2, respectively.
The electron and ion thermal velocities are determined from the pressure balance, $(v_{t,e},v_{t,i}) = (0.1c/\sqrt{2},0.02c/\sqrt{2})$, and the current sheet thickness is chosen to be equal to the ion inertia length.
The simulation domain in physical space spanning $[0,2560\lambda_D] \times [0,640\lambda_D]$ is discretized by $256\times64$ grid points, where the Debye length $\lambda_D$ is defined in the current sheet.
The time step is $0.25\omega_{pe}^{-1}$.

\begin{figure}[htbp]
 \centering
\includegraphics[clip,angle=0,scale=0.6]{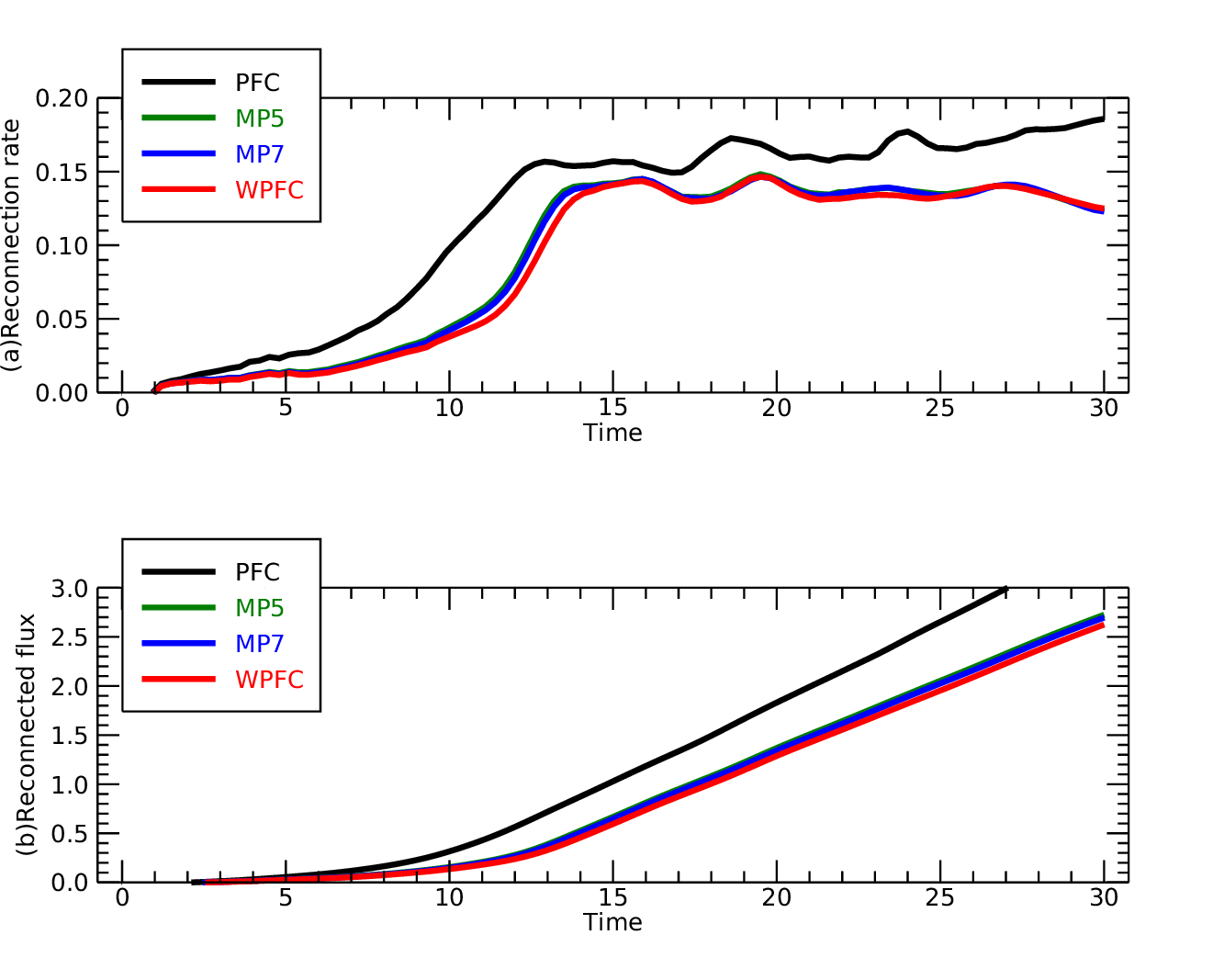}
\caption{Time profiles of (a) the reconnection rate and (b) the reconnected flux obtained from the third-order PFC scheme (black), the fifth- and seventh-order monotonicity-preserving schemes (green and blue), and the present scheme (red), respectively. The blue lines overlap with the green lines.}
\label{fig:mrate}
\end{figure}

\begin{figure}[htbp]
 \centering
\includegraphics[clip,angle=0,scale=0.5]{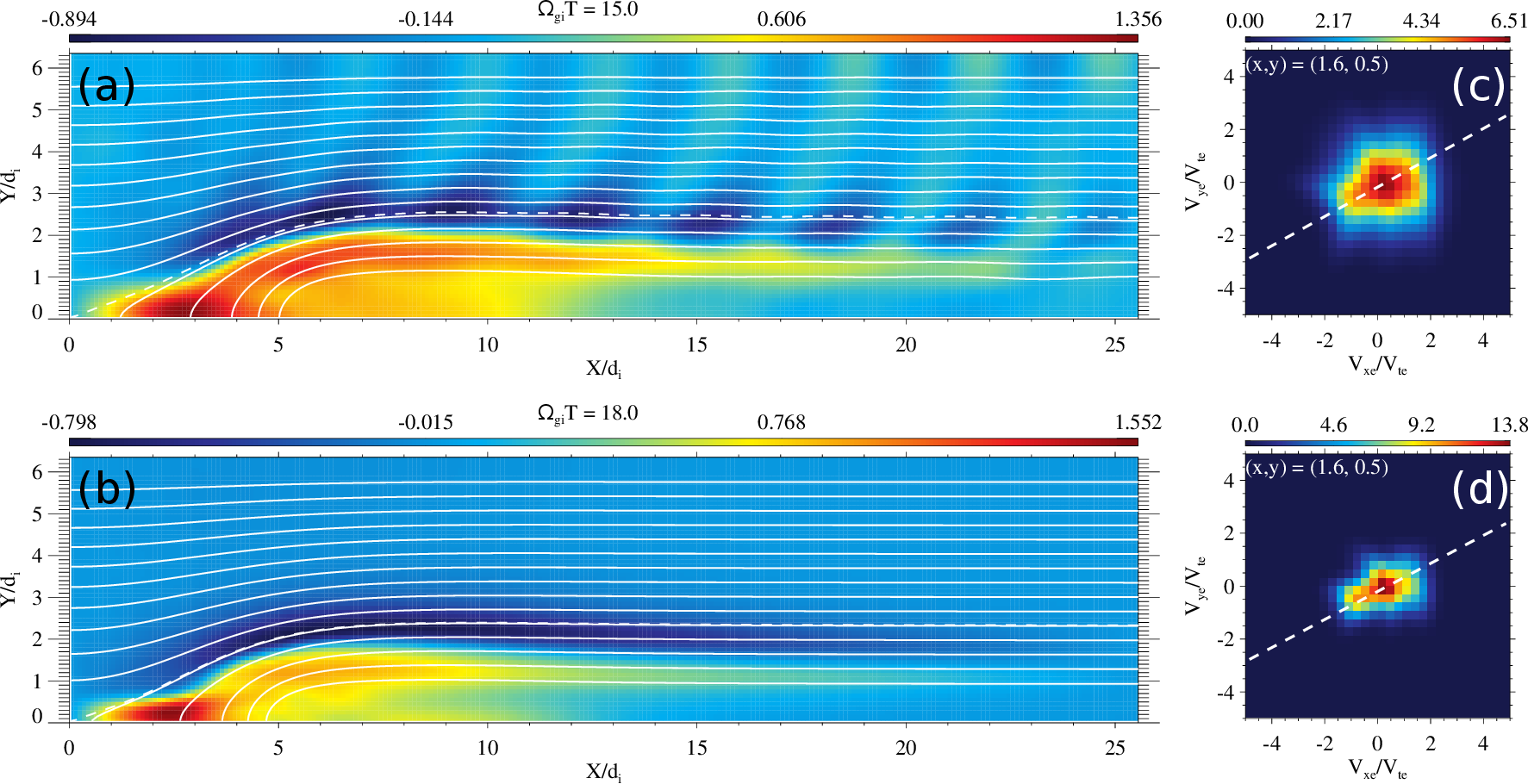}
\caption{Numerical experiments of the magnetic reconnection obtained from (top) the third-order PFC and (bottom) the present schemes. (a,b) Spatial distribution of the $x$-component of the electron bulk velocity normalized by the upstream {\Alfven} velocity. Solid and dashed lines indicate the magnetic field line and separatrix, respectively. (c,d) Electron velocity distribution function $f_e(v_x,v_y)$ at $(x,y)=(1.6,0.5)$. Dashed lines indicate the magnetic field direction.}
\label{fig:mrx}
\end{figure}

The reconnection is triggered by imposing a localized magnetic field perturbation around $(x,y)=(0,0)$.
Following the approach by \cite{2011PhPl...18l2108Z}, we identify the dominant reconnection site, and compute the normalized reconnection rate \reviser{and the amount of reconnected flux by measuring the magnetic flux there.}
Figure \ref{fig:mrate}(a) shows the time profile of the reconnection rate.
 The result obtained from the PFC scheme exhibits earlier growth by $t=13\Omega_{gi}^{-1}$ where $\Omega_{gi}$ is the ion gyrofrequency, and continues to increase thereafter.
In contrast, the rates obtained from the high-order schemes saturate at $t=15\Omega_{gi}^{-1}$ with a magnitude of approximately 0.14, consistent with the previous study \reviser{by \cite{2011PhPl...18l2108Z}}.
Figure \ref{fig:mrx} compares the $x$-component of the electron bulk velocity and the electron velocity distribution function obtained from the PFC scheme at $t=15\Omega_{gi}^{-1}$ (top panels) and the present scheme at $t=18\Omega_{gi}^{-1}$ (bottom panels).
\reviser{Because the reconnection rates differ between the two schemes, the comparison is made at the time when the reconnected flux reaches 1.0 (Figure \ref{fig:mrate}(b)).}
In Figure \ref{fig:mrx}(a), the PFC scheme exhibits a broader spatial distribution in the $y$-direction and a wave structure along the $x$-direction in the upstream region, compared to the present scheme in Figure \ref{fig:mrx}(b) and the MP schemes (not shown).
The electron velocity distribution near the reconnection site $(x,y)=(1.6,0.5)$ is elongated primarily along the magnetic field line due to the decrease in magnetic field strength.
However, considerable broadening perpendicular to the magnetic field line is observed in the PFC scheme (Figure \ref{fig:mrx}(c)), attributed to numerical diffusion when solving the gyromotion.
This artificial perpendicular electron heating starts in the upstream region, and leads to the unphysical generation of plasma waves propagating along the magnetic field line as seen in Figure \ref{fig:mrx}(a).

\section{Conclusion} \label{sec:conclusion}
We presented a weighted positive and flux conservative (WPFC) method for solving the Vlasov equation, which builds upon the third-order PFC method developed by \cite{2001JCoPh.172..166F}, extending it to attain fifth-order accuracy.
This is achieved by a convex combination of positive and non-oscillatory polynomials in substencils, inspired by the weighted essentially non-oscillatory (WENO) scheme.
Unlike conventional WENO schemes, however, our scheme formulates nonlinear weights for these polynomials to assign higher priority to substencils with larger $L^2$ norm, that is, steeper substencils.
This approach enhances resolution while maintaining positivity and non-oscillatory properties.
An approximate dissipation relation reveals that the spectral properties of the present scheme outperform those of the underlying fifth-order scheme and even surpass those of a seventh-order scheme in certain wavenumber ranges.
Numerical experiments on linear advection suggest that well-designed nonlinear weights can enhance both the accuracy and stability of the scheme, outperforming the underlying linear scheme.

We applied the present scheme to the electrostatic Vlasov-\Ampere simulations and demonstrated improved conservation of entropy.
In the simulations of the two stream instability and the bump-on-tail instability, we show that the high-order scheme significantly improves conservation properties compared to the low-order scheme, supporting the advantage of high-order schemes.
The present scheme exhibits better conservation of entropy compared to the fifth- and seventh-order monotonicity-preserving schemes.
In the simulation of the nonlinear Landau damping, in contrast, the conservation of entropy is almost the same level among the low- and high-order schemes, with the low-order scheme providing the best energy conservation.
This implies that utilizing higher-order schemes does not always improve conservation properties in certain problems, possibly due to more pronounced occurrence and subsequent numerical dissipation of grid-scale structure in velocity space.
Our approach of nonlinear weighting based on the $L^2$ norm effectively improves the conservation in the Vlasov-\Ampere simulations, even though the formal order of accuracy remains unchanged.

\begin{figure}[htbp]
\centering
\includegraphics[clip,angle=0,scale=0.7]{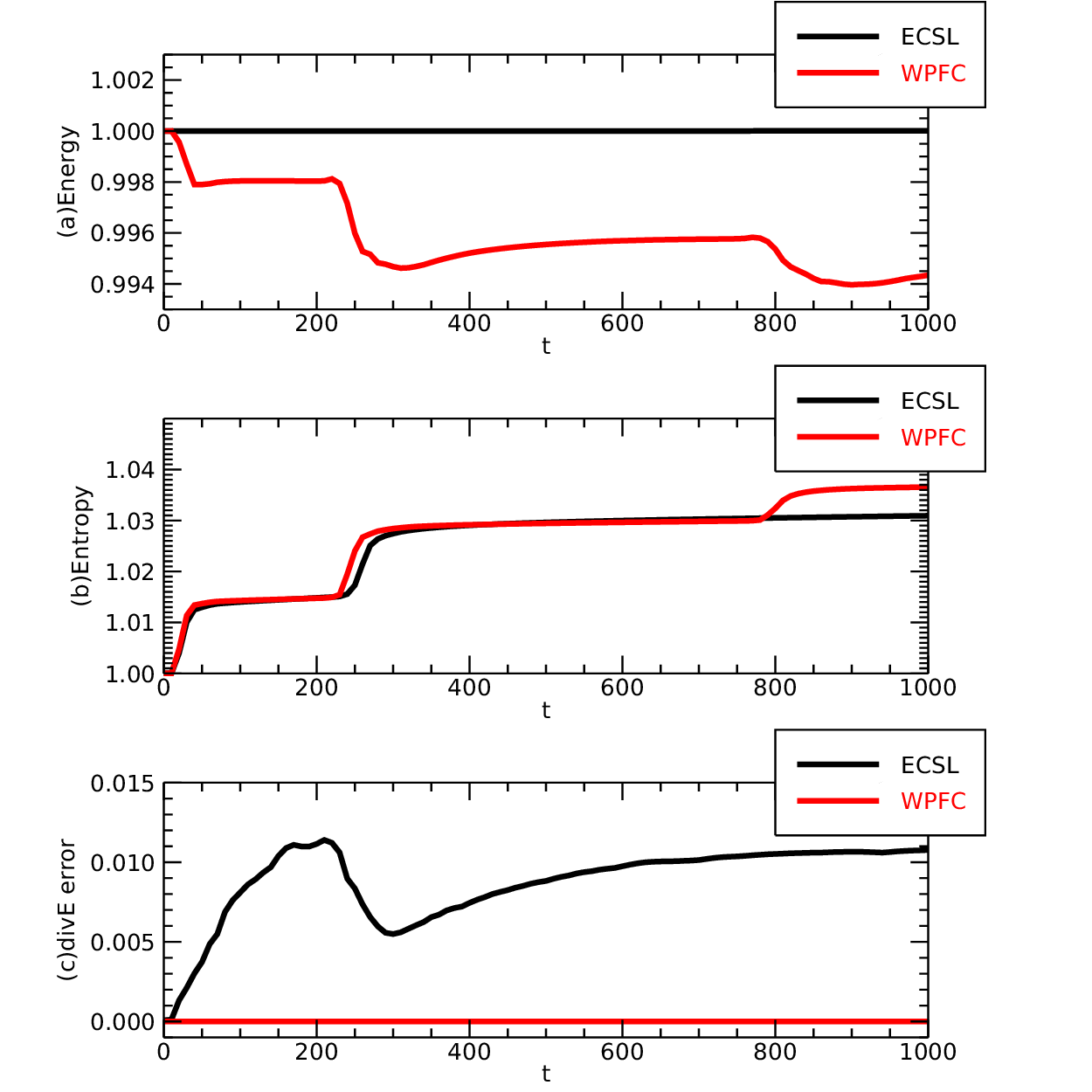} 
\caption{\revrr{Comparison between the energy conserving semi-Lagrangian (ECSL) scheme (black lines) and the present scheme (red lines) in the two-stream instability problem: (a) total energy, (b) entropy, and (c) absolute error in Gauss's law for the electric field. The increase in entropy at $t=800$ in the present scheme corresponds to the transition to the single phase space hole (Figure \ref{fig:twost}(b)), which did not occurred in the ECSL run.}}
\label{fig:ecsl_wpfc}
\end{figure}

\reviseb{In terms of energy conservation, the energy conserving semi-Lagrangian (ECSL) scheme \citep{liu2023efficient} provides improved accuracy.
The scheme introduces a novel electric field update by coupling the \Ampere and Vlasov momentum equations in an efficient semi-implicit manner, thereby ensuring total energy conservation at the discrete level.
\revrr{Moreover, in contrast to explicit schemes including the present one, the ECSL scheme provides reliable solutions even when the time step does not fully resolve the plasma frequency, which is advantageous for long-term simulations.
Meanwhile, the scheme does not guarantee a non-oscillatory property, because nonlinear schemes to solve the advection equation in velocity space without violating energy conservation are not yet available.}
\revrr{In Figure \ref{fig:ecsl_wpfc}, we applied the fifth-order ECSL scheme to the two stream instability problem and confirmed that total energy is conserved, while entropy conservation is comparable to that of the present scheme and Gauss's law for the electric field is not strictly satisfied.}
\revrr{These advantages and drawbacks represent a trade-off, and simultaneously ensuring numerical stability constraints and strict physical constraints in long-term simulations remains a challenge.}}

We also applied the scheme to the Vlasov-Maxwell simulation of collisionless magnetic reconnection.
The solution is considerably affected by the order of accuracy used in velocity space.
In particular, the low-order scheme leads to artificial plasma heating, resulting in the degradation of the global evolution and the spatial distribution, even when the resolution in physical space is fixed.
This result supports the importance of employing high-order schemes in velocity space even for multidimensional problems, enhancing the significance of the present study.

\reviser{A drawback of the present scheme is the computational cost arising from the calculation of the nonlinear reconstruction functions in three substencils and their nonlinear convex combination.
To assess computational efficiency, we measured the elapsed time for the two stream instability problem.
The simulations were performed on dual Intel Xeon Gold 6326 (2.9 GHz) processors, with the code compiled using the GNU C++ Compiler and optimization flags \texttt{-march=native -flto -fopenmp}, together with either \texttt{-O2} or \texttt{-O3}.
With \texttt{-O2}, the measured elapsed times of the PFC, MP5, MP7, and WPFC schemes were 21.9, 21.5, 22.1, and 25.2 seconds, respectively (ratios of $1.02 : 1.0 : 1.03 : 1.17$). With \texttt{-O3}, the results changed to 23.1, 20.6, 20.8, and 21.7 seconds ($1.12 : 1.0 : 1.01 : 1.05$).
High-level optimization improves the computational efficiency of the high-order schemes, especially the present scheme.
Further aggressive optimization will improve the practicality of the present scheme.}

Given the similarity of our approach to the improved WENO-Z schemes, we expect that it is not limited solely to the Vlasov simulation but is useful to other hyperbolic systems as well.
The application of the present scheme to hydrodynamic simulations is planned for future research.

\section*{Declarations}

\section*{Availability of data and materials}
The datasets used and/or analyzed during the current study are available from the corresponding author upon reasonable request.

\section*{Competing interests}
The authors declare that they have no competing interests.

\section*{Funding}
JSPS KAKENHI Grant No. JP20K04056 and JP25H00625.
 
\section*{Authors' contributions}
TM developed the concept and methodology of this study, and performed numerical simulations with support from YM. 
All authors contributed to the interpretation of the results.
TM wrote the original draft, and YM reviewed and edited it.


\acknowledgments{
We thank the anonymous referees for carefully reviewing our manuscript and providing valuable comments.
}

\bibliographystyle{plainnat}

\begin{thebibliography}{36}
\providecommand{\natexlab}[1]{#1}
\providecommand{\url}[1]{\texttt{#1}}
\expandafter\ifx\csname urlstyle\endcsname\relax
  \providecommand{\doi}[1]{doi: #1}\else
  \providecommand{\doi}{doi: \begingroup \urlstyle{rm}\Url}\fi

\bibitem[Acker et~al.(2016)Acker, Borges, and Costa]{acker2016improved}
Felipe Acker, RB~de~R Borges, and Bruno Costa.
\newblock {An improved WENO-Z scheme}.
\newblock \emph{Journal of Computational Physics}, 313:\penalty0 726--753,
  2016.

\bibitem[{Birdsall} and {Langdon}(1991)]{PIC}
C.~K. {Birdsall} and A.~B. {Langdon}.
\newblock \emph{{Plasma Physics via Computer Simulation}}.
\newblock Inst. of Phys. Publishing, Bristol/Philadelphia, 1991.

\bibitem[{Birn} and {Hesse}(2001)]{2001JGR...106.3737B}
J.~{Birn} and M.~{Hesse}.
\newblock {Geospace Environment Modeling (GEM) magnetic reconnection challenge:
  Resistive tearing, anisotropic pressure and hall effects}.
\newblock \emph{\jgr}, 106:\penalty0 3737--3750, March 2001.
\newblock \doi{10.1029/1999JA001001}.

\bibitem[{Borges} et~al.(2008){Borges}, {Carmona}, {Costa}, and
  {Don}]{2008JCoPh.227.3191B}
R.~{Borges}, M.~{Carmona}, B.~{Costa}, and W.~S. {Don}.
\newblock {An improved weighted essentially non-oscillatory scheme for
  hyperbolic conservation laws}.
\newblock \emph{Journal of Computational Physics}, 227:\penalty0 3191--3211,
  March 2008.
\newblock \doi{10.1016/j.jcp.2007.11.038}.

\bibitem[Cai et~al.(2016)Cai, Qiu, and Qiu]{cai2016conservative}
Xiaofeng Cai, Jianxian Qiu, and Jing-Mei Qiu.
\newblock {A conservative semi-Lagrangian HWENO method for the Vlasov
  equation}.
\newblock \emph{Journal of Computational Physics}, 323:\penalty0 95--114, 2016.

\bibitem[{Cheng} and {Knorr}(1976)]{1976JCoPh..22..330C}
C.~Z. {Cheng} and G.~{Knorr}.
\newblock {The integration of the Vlasov equation in configuration space}.
\newblock \emph{Journal of Computational Physics}, 22:\penalty0 330--351,
  November 1976.
\newblock \doi{10.1016/0021-9991(76)90053-X}.

\bibitem[{Crouseilles} et~al.(2009){Crouseilles}, {Respaud}, and
  {Sonnendr{\"u}cker}]{2009CoPhC.180.1730C}
N.~{Crouseilles}, T.~{Respaud}, and E.~{Sonnendr{\"u}cker}.
\newblock {A forward semi-Lagrangian method for the numerical solution of the
  Vlasov equation}.
\newblock \emph{Computer Physics Communications}, 180:\penalty0 1730--1745,
  October 2009.
\newblock \doi{10.1016/j.cpc.2009.04.024}.

\bibitem[{Crouseilles} et~al.(2010){Crouseilles}, {Mehrenberger}, and
  {Sonnendr{\"u}cker}]{2010JCoPh.229.1927C}
N.~{Crouseilles}, M.~{Mehrenberger}, and E.~{Sonnendr{\"u}cker}.
\newblock {Conservative semi-Lagrangian schemes for Vlasov equations}.
\newblock \emph{Journal of Computational Physics}, 229:\penalty0 1927--1953,
  March 2010.
\newblock \doi{10.1016/j.jcp.2009.11.007}.

\bibitem[{Filbet} and {Sonnendr{\"u}cker}(2003)]{2003CoPhC.150..247F}
F.~{Filbet} and E.~{Sonnendr{\"u}cker}.
\newblock {Comparison of Eulerian Vlasov solvers}.
\newblock \emph{Computer Physics Communications}, 150:\penalty0 247--266,
  February 2003.
\newblock \doi{10.1016/S0010-4655(02)00694-X}.

\bibitem[{Filbet} et~al.(2001){Filbet}, {Sonnendr{\"u}cker}, and
  {Bertrand}]{2001JCoPh.172..166F}
F.~{Filbet}, E.~{Sonnendr{\"u}cker}, and P.~{Bertrand}.
\newblock {Conservative Numerical Schemes for the Vlasov Equation}.
\newblock \emph{Journal of Computational Physics}, 172:\penalty0 166--187,
  September 2001.
\newblock \doi{10.1006/jcph.2001.6818}.

\bibitem[{Jiang} and {Shu}(1996)]{1996JCoPh.126..202J}
G.~{Jiang} and C.~{Shu}.
\newblock {Efficient Implementation of Weighted ENO Schemes}.
\newblock \emph{Journal of Computational Physics}, 126:\penalty0 202--228, June
  1996.
\newblock \doi{10.1006/jcph.1996.0130}.

\bibitem[{Lele}(1992)]{1992JCoPh.103...16L}
S.~K. {Lele}.
\newblock {Compact Finite Difference Schemes with Spectral-like Resolution}.
\newblock \emph{Journal of Computational Physics}, 103:\penalty0 16--42,
  November 1992.
\newblock \doi{10.1016/0021-9991(92)90324-R}.

\bibitem[Liu et~al.(2023)Liu, Cai, Cao, and Lapenta]{liu2023efficient}
Hongtao Liu, Xiaofeng Cai, Yong Cao, and Giovanni Lapenta.
\newblock {An efficient energy conserving semi-Lagrangian kinetic scheme for
  the Vlasov-Amp{\`e}re system}.
\newblock \emph{Journal of Computational Physics}, 492:\penalty0 112412, 2023.

\bibitem[{Liu} et~al.(2025){Liu}, {Lu}, {Xia}, {Keppens}, and
  {Lapenta}]{2025JCoPh.52913858L}
Hongtao {Liu}, Chang {Lu}, Guangqing {Xia}, Rony {Keppens}, and Giovanni
  {Lapenta}.
\newblock {An efficient energy conserving semi-Lagrangian kinetic scheme for
  the Vlasov-Maxwell system}.
\newblock \emph{Journal of Computational Physics}, 529:\penalty0 113858, May
  2025.
\newblock \doi{10.1016/j.jcp.2025.113858}.

\bibitem[Luo and Wu(2021)]{luo2021improved}
Xin Luo and Song-ping Wu.
\newblock {An improved WENO-Z+ scheme for solving hyperbolic conservation
  laws}.
\newblock \emph{Journal of Computational Physics}, 445:\penalty0 110608, 2021.

\bibitem[{Minoshima} et~al.(2011){Minoshima}, {Matsumoto}, and
  {Amano}]{2011JCoPh.230.6800M}
T.~{Minoshima}, Y.~{Matsumoto}, and T.~{Amano}.
\newblock {Multi-moment advection scheme for Vlasov simulations}.
\newblock \emph{Journal of Computational Physics}, 230:\penalty0 6800--6823,
  July 2011.
\newblock \doi{10.1016/j.jcp.2011.05.010}.

\bibitem[{Minoshima} et~al.(2013){Minoshima}, {Matsumoto}, and
  {Amano}]{2013JCoPh.236...81M}
T.~{Minoshima}, Y.~{Matsumoto}, and T.~{Amano}.
\newblock {Multi-moment advection scheme in three dimension for Vlasov
  simulations of magnetized plasma}.
\newblock \emph{Journal of Computational Physics}, 236:\penalty0 81--95, March
  2013.
\newblock \doi{10.1016/j.jcp.2012.11.024}.

\bibitem[{Minoshima} et~al.(2015){Minoshima}, {Matsumoto}, and
  {Amano}]{2015CoPhC.187..137M}
T.~{Minoshima}, Y.~{Matsumoto}, and T.~{Amano}.
\newblock {A finite volume formulation of the multi-moment advection scheme for
  Vlasov simulations of magnetized plasma}.
\newblock \emph{Computer Physics Communications}, 187:\penalty0 137--151,
  February 2015.
\newblock \doi{10.1016/j.cpc.2014.10.023}.

\bibitem[{Minoshima} et~al.(2019){Minoshima}, {Miyoshi}, and
  {Matsumoto}]{2019ApJS..242...14M}
Takashi {Minoshima}, Takahiro {Miyoshi}, and Yosuke {Matsumoto}.
\newblock {A High-order Weighted Finite Difference Scheme with a Multistate
  Approximate Riemann Solver for Divergence-free Magnetohydrodynamic
  Simulations}.
\newblock \emph{\apjs}, 242\penalty0 (2):\penalty0 14, Jun 2019.
\newblock \doi{10.3847/1538-4365/ab1a36}.

\bibitem[{Nakamura} and {Yabe}(1999)]{1999CoPhC.120..122N}
T.~{Nakamura} and T.~{Yabe}.
\newblock {Cubic interpolated propagation scheme for solving the
  hyper-dimensional Vlasov-Poisson equation in phase space}.
\newblock \emph{Computer Physics Communications}, 120:\penalty0 122--154,
  August 1999.
\newblock \doi{10.1016/S0010-4655(99)00247-7}.

\bibitem[Pirozzoli(2006)]{pirozzoli2006spectral}
Sergio Pirozzoli.
\newblock {On the spectral properties of shock-capturing schemes}.
\newblock \emph{Journal of Computational Physics}, 219\penalty0 (2):\penalty0
  489--497, 2006.

\bibitem[Qiu and Christlieb(2010)]{qiu2010conservative}
Jing-Mei Qiu and Andrew Christlieb.
\newblock {A conservative high order semi-Lagrangian WENO method for the Vlasov
  equation}.
\newblock \emph{Journal of Computational Physics}, 229\penalty0 (4):\penalty0
  1130--1149, 2010.

\bibitem[{Schmitz} and {Grauer}(2006{\natexlab{a}})]{2006CoPhC.175...86S}
H.~{Schmitz} and R.~{Grauer}.
\newblock {Comparison of time splitting and backsubstitution methods for
  integrating Vlasov's equation with magnetic fields}.
\newblock \emph{Computer Physics Communications}, 175:\penalty0 86--92, July
  2006{\natexlab{a}}.
\newblock \doi{10.1016/j.cpc.2006.02.007}.

\bibitem[{Schmitz} and {Grauer}(2006{\natexlab{b}})]{2006PhPl...13i2309S}
H.~{Schmitz} and R.~{Grauer}.
\newblock {Kinetic Vlasov simulations of collisionless magnetic reconnection}.
\newblock \emph{Physics of Plasmas}, 13\penalty0 (9):\penalty0 092309--+,
  September 2006{\natexlab{b}}.
\newblock \doi{10.1063/1.2347101}.

\bibitem[{Shoucri}(2017)]{2017LPB....35..706S}
Magdi {Shoucri}.
\newblock {Formation of electron holes in the long-time evolution of the
  bump-on-tail instability}.
\newblock \emph{Laser and Particle Beams}, 35\penalty0 (4):\penalty0 706--721,
  December 2017.
\newblock \doi{10.1017/S0263034617000775}.

\bibitem[Sirajuddin and Hitchon(2019)]{sirajuddin2019truly}
David Sirajuddin and William~NG Hitchon.
\newblock {A truly forward semi-Lagrangian WENO scheme for the Vlasov-Poisson
  system}.
\newblock \emph{Journal of Computational Physics}, 392:\penalty0 619--665,
  2019.

\bibitem[{Sircombe} and {Arber}(2009)]{2009JCoPh.228.4773S}
N.~J. {Sircombe} and T.~D. {Arber}.
\newblock {VALIS: A split-conservative scheme for the relativistic 2D
  Vlasov-Maxwell system}.
\newblock \emph{Journal of Computational Physics}, 228:\penalty0 4773--4788,
  July 2009.
\newblock \doi{10.1016/j.jcp.2009.03.029}.

\bibitem[{Sonnendr{\"u}cker} et~al.(1999){Sonnendr{\"u}cker}, {Roche},
  {Bertrand}, and {Ghizzo}]{1999JCoPh.149..201S}
E.~{Sonnendr{\"u}cker}, J.~{Roche}, P.~{Bertrand}, and A.~{Ghizzo}.
\newblock {The Semi-Lagrangian Method for the Numerical Resolution of the
  Vlasov Equation}.
\newblock \emph{Journal of Computational Physics}, 149:\penalty0 201--220,
  March 1999.
\newblock \doi{10.1006/jcph.1998.6148}.

\bibitem[{Suresh} and {Huynh}(1997)]{1997JCoPh.136...83S}
A.~{Suresh} and H.~T. {Huynh}.
\newblock {Accurate Monotonicity-Preserving Schemes with Runge Kutta Time
  Stepping}.
\newblock \emph{Journal of Computational Physics}, 136:\penalty0 83--99,
  September 1997.
\newblock \doi{10.1006/jcph.1997.5745}.

\bibitem[{Tanaka} et~al.(2017){Tanaka}, {Yoshikawa}, {Minoshima}, and
  {Yoshida}]{2017ApJ...849...76T}
Satoshi {Tanaka}, Kohji {Yoshikawa}, Takashi {Minoshima}, and Naoki {Yoshida}.
\newblock {Multidimensional Vlasov-Poisson Simulations with High-order
  Monotonicity- and Positivity-preserving Schemes}.
\newblock \emph{\apj}, 849\penalty0 (2):\penalty0 76, Nov 2017.
\newblock \doi{10.3847/1538-4357/aa901f}.

\bibitem[{Tsallis}(1988)]{1988JSP....52..479T}
C.~{Tsallis}.
\newblock {Possible generalization of Boltzmann-Gibbs statistics}.
\newblock \emph{Journal of Statistical Physics}, 52:\penalty0 479--487, July
  1988.
\newblock \doi{10.1007/BF01016429}.

\bibitem[{Umeda}(2008)]{2008EP&S...60..773U}
T.~{Umeda}.
\newblock {A conservative and non-oscillatory scheme for Vlasov code
  simulations}.
\newblock \emph{Earth, Planets, and Space}, 60:\penalty0 773--779, July 2008.

\bibitem[{Umeda} et~al.(2009){Umeda}, {Togano}, and
  {Ogino}]{2009CoPhC.180..365U}
T.~{Umeda}, K.~{Togano}, and T.~{Ogino}.
\newblock {Two-dimensional full-electromagnetic Vlasov code with conservative
  scheme and its application to magnetic reconnection}.
\newblock \emph{Computer Physics Communications}, 180:\penalty0 365--374, March
  2009.
\newblock \doi{10.1016/j.cpc.2008.11.001}.

\bibitem[{Umeda} et~al.(2012){Umeda}, {Nariyuki}, and
  {Kariya}]{2012CoPhC.183.1094U}
Takayuki {Umeda}, Yasuhiro {Nariyuki}, and Daichi {Kariya}.
\newblock {A non-oscillatory and conservative semi-Lagrangian scheme with
  fourth-degree polynomial interpolation for solving the Vlasov equation}.
\newblock \emph{Computer Physics Communications}, 183\penalty0 (5):\penalty0
  1094--1100, May 2012.
\newblock \doi{10.1016/j.cpc.2012.01.011}.

\bibitem[Xiong et~al.(2014)Xiong, Qiu, Xu, and Christlieb]{xiong2014high}
Tao Xiong, Jing-Mei Qiu, Zhengfu Xu, and Andrew Christlieb.
\newblock {High order maximum principle preserving semi-Lagrangian finite
  difference WENO schemes for the Vlasov equation}.
\newblock \emph{Journal of Computational Physics}, 273:\penalty0 618--639,
  2014.

\bibitem[{Zenitani} et~al.(2011){Zenitani}, {Hesse}, {Klimas}, {Black}, and
  {Kuznetsova}]{2011PhPl...18l2108Z}
S.~{Zenitani}, M.~{Hesse}, A.~{Klimas}, C.~{Black}, and M.~{Kuznetsova}.
\newblock {The inner structure of collisionless magnetic reconnection: The
  electron-frame dissipation measure and Hall fields}.
\newblock \emph{Physics of Plasmas}, 18\penalty0 (12):\penalty0 122108,
  December 2011.
\newblock \doi{10.1063/1.3662430}.

\end{thebibliography}

\clearpage

\clearpage

\clearpage

\clearpage

\clearpage

\clearpage

\clearpage

\clearpage

\clearpage

\clearpage
\begin{table}
\centering
\caption{Accuracy on the linear advection of the Gaussian profile.} 
\begin{tabular}{cccccc}\hline \hline
 Scheme& Grids & $L^1$ error & $L^1$ order & $L^{\infty}$ error & $L^{\infty}$ order\\ \hline
Linear & 32 & 4.79e-2 & $-$ & 2.56e-1 & $-$ \\
 & 64 & 1.01e-2 & 2.25 & 8.49e-2 & 1.59 \\
 & 128 & 6.88e-4 & 3.88 & 8.39e-3 & 3.34 \\
 & 256 & 2.42e-5 & 4.83 & 3.29e-4 & 4.67 \\ \hline
MP5 & 32 & 3.60e-2 & $-$ & 2.92e-1 & $-$ \\
 & 64 & 6.81e-3 & 2.40 & 8.03e-2 & 1.86 \\
 & 128 & 6.88e-4 & 3.31 & 8.39e-3 & 3.26 \\
 & 256 & 2.42e-5 & 4.83 & 3.29e-4 & 4.67 \\ \hline
WPFC & 32 & 3.27e-2 & $-$ & 2.63e-1 & $-$ \\
 & 64 & 4.47e-3 & 2.87 & 4.23e-2 & 2.63 \\
 & 128 & 4.29e-4 & 3.38 & 3.72e-3 & 3.51 \\
 & 256 & 1.47e-5 & 4.87 & 1.31e-4 & 4.83 \\ \hline
\end{tabular}
\label{tab:acc1}
\end{table}

\begin{table}
\centering
\caption{Accuracy on the linear advection of the sinusoidal profile.} 
\begin{tabular}{cccccc}\hline \hline
 Scheme& Grids & $L^1$ error & $L^1$ order & $L^{\infty}$ error & $L^{\infty}$ order\\ \hline
Linear & 32 & 2.39e-2 & $-$ & 3.47e-2 & $-$ \\
 & 64 & 8.52e-4 & 4.81 & 1.31e-3 & 4.73 \\
 & 128 & 2.71e-5 & 4.97 & 4.23e-5 & 4.95 \\
 & 256 & 8.49e-7 & 5.00 & 1.33e-6 & 4.99 \\ \hline
WPFC & 32 & 1.38e-3 & $-$ & 1.83e-3 & $-$ \\
 & 64 & 8.17e-5 & 4.08 & 1.24e-4 & 3.88 \\
 & 128 & 2.06e-6 & 5.31 & 4.99e-6 & 4.64 \\
 & 256 & 2.83e-8 & 6.19 & 1.34e-7 & 5.22 \\ \hline
\end{tabular}
\label{tab:acc2}
\end{table}

\begin{table}
 \centering
\caption{Accuracy on the solid body rotation of the Gaussian profile.}
\begin{tabular}{cccccc}\hline \hline
 Scheme& Grids & $L^1$ error & $L^1$ order & $L^{\infty}$ error & $L^{\infty}$ order\\ \hline
Linear & $32^3$ & 9.03e-4 & $-$ & 9.03e-2 & $-$ \\ 
& $64^3$ & 4.90e-5 & 4.20 & 6.91e-3 & 3.71 \\ 
& $128^3$ & 1.79e-6 & 4.77 & 2.49e-4 & 4.79 \\ 
\hline
MP5 & $32^3$ & 6.71e-4 & $-$ & 8.95e-2 & $-$ \\ 
& $64^3$ & 4.90e-5 & 3.78 & 6.91e-3 & 3.70 \\ 
& $128^3$ & 1.79e-6 & 4.77 & 2.49e-4 & 4.79 \\ 
\hline
WPFC & $32^3$ & 4.89e-4 & $-$ & 4.46e-2 & $-$ \\ 
& $64^3$ & 3.84e-5 & 3.67 & 2.88e-3 & 3.95 \\ 
& $128^3$ & 1.40e-6 & 4.78 & 1.12e-4 & 4.68 \\ 
\hline
\end{tabular}
\label{tab:acc3}
\end{table}

\begin{table}
 \centering
\caption{Same as Table \ref{tab:acc3}, but using smaller time steps $\Delta t \propto (\Delta x)^{2.5}$.}
\begin{tabular}{cccccc}\hline \hline
 Scheme& Grids & $L^1$ error & $L^1$ order & $L^{\infty}$ error & $L^{\infty}$ order\\ \hline
WPFC & $32^3$ & 4.89e-4 & $-$ & 4.46e-2 & $-$ \\ 
& $64^3$ & 3.86e-5 & 3.66 & 2.88e-3 & 3.95 \\ 
& $128^3$ & 1.30e-6 & 4.89 & 9.66e-5 & 4.89 \\ 
\hline
\end{tabular}
\label{tab:acc4}
\end{table}

\end{document}